\documentclass{article}
\usepackage{amsmath,amsthm}
\usepackage{verbatim,amsmath,amscd,amssymb,amsthm}
\theoremstyle{definition}
\newtheorem{df}{Definition}[section]
\newtheorem{rem}[df]{Remark}

\theoremstyle{plain}

\newtheorem{pro}[df]{Proposition}
\newtheorem{lem}[df]{Lemma}
\newtheorem{thm}[df]{Theorem}
\newtheorem{cor}[df]{Corollary}
\begin{document}
\font\germ=eufm10
\def\ssl{\hbox{\germ sl}}
\def\slh{\widehat{\ssl_2}}
\def\ge{\hbox{\germ g}}
\def\aaa{@}
\title{\Large\bf Tensor Product of the Fundamental Representations 
for the Quantum Loop Algebras of Type $A$ at Roots of Unity}

\author{
Yuuki A\textsc{be} 
\thanks
{
e-mail: abeyuuki@ms.u-tokyo.ac.jp
}
\\
Graduate School of Mathematical Sciences, 
\\
University of Tokyo 
}
\date{}
\maketitle
\begin{abstract}
In this paper, 
we consider the necessary and sufficient conditions 
for the tensor product of the fundamental representations 
for the restricted quantum loop algebras of type $A$ 
at roots of unity to be irreducible.
\end{abstract}
\maketitle

\renewcommand{\labelenumi}{$($\roman{enumi}$)$}
\renewcommand{\labelenumii}{$(${\rm \alph{enumii}}$)$}
\font\germ=eufm10
\newcommand{\fra}{\mathfrak a}
\newcommand{\frg}{\mathfrak g}
\def\gl{{\mathfrak g}{\mathfrak l}}
\newcommand{\frh}{\mathfrak h}
\newcommand{\frp}{\mathfrak p}
\newcommand{\frs}{\mathfrak s}
\def\sl{{\mathfrak s}{\mathfrak l}}
\def\h{\frh}
\newcommand{\bbC}{\mathbb C}
\newcommand{\bbN}{\mathbb N}
\newcommand{\bbQ}{\mathbb Q}
\newcommand{\bbR}{\mathbb R}
\newcommand{\bbZ}{\mathbb Z}
\def\AA{{\cal A}}
\def\DD{{\cal D}}
\def\JJ{{\cal J}}
\def\PP{{\cal P}}
\def\RR{{\cal R}}
\def\SS{{\cal S}}
\def\TT{{\cal T}}
\def\WW{{\cal W}}
\def\al{\alpha}
\def\be{\beta}
\def\de{\delta}
\def\De{\Delta}
\def\ep{\epsilon}
\def\ga{\gamma}
\def\la{\lambda}
\def\La{\Lambda}
\def\om{\omega}
\def\Om{\Omega}
\def\th{\theta}
\def\vep{\varepsilon}
\def\vph{\varphi}
\def\ze{\zeta}
\def\bf{\textbf}
\def\it{\textit}
\def\mbf{\mathbf}
\def\mrm{\mathrm}
\def\rm{\textrm}
\def\arr{\longrightarrow}
\def\bs{\backslash}
\def\no{\nonumber}
\def\ol{\overline}
\def\ot{\otimes}
\def\q{\quad}
\def\qq{\qquad}
\def\s={\star \in \{-,0,+\}}
\def\wh{\widehat}
\def\wt{\widetilde}
\def\BA{B_{\AA}}
\def\Bfin{B_{\vep}^{\mrm{fin}}}
\def\Bres{B_{\vep}^{\mrm{res}}}
\def\EUAres{U_{\AA}^{\mrm{res'}}}
\def\EUfin{U_{\vep}^{\mrm{fin'}}}
\def\EUres{U_{\vep}^{\mrm{res'}}}
\def\Lfin{L_{\vep}^{\mrm{fin}}}
\def\Lnil{L_{\vep}^{\mrm{nil}}}
\def\Pfin{P_{\vep}^{\mrm{fin}}}
\def\Pnil{P_{\vep}^{\mrm{nil}}}
\def\Pres{P_{\vep}^{\mrm{res}}}
\def\slEAres{\EUAres(\sl_{n+1})}
\def\slEfin{\EUfin(\sl_{n+1})}
\def\slEres{\EUres(\sl_{n+1})}
\def\slEU{\Ue^{'}(\sl_{n+1})}
\def\slfin{\Ufin(\sl_{n+1})}
\def\slres{\Ures(\sl_{n+1})}
\def\slU{\Ue(\sl_{n+1})}
\def\slUA{\UA(\sl_{n+1})}
\def\slUAres{\UAres(\sl_{n+1})}
\def\tBA{\wt{B}_{\cal{A}}}
\def\tBAres{\wt{B}_{{\cal A}}^{\mrm{res}}}
\def\tBres{\wt{B}_{\vep}^{\mrm{res}}}
\def\tBfin{\tB_{\vep}^{\mrm{fin}}}
\def\tPres{\wt{P}_{\vep}^{\mrm{res}}}
\def\tPfin{\wt{P}_{\vep}^{\mrm{fin}}}
\def\tslfin{\Ufin(\wt{\sl}_{n+1})}
\def\tslres{\Ures(\wt{\sl}_{n+1})}
\def\tslU{\Ue(\wt{\sl}_{n+1})}
\def\tUA{\wt{U}_{\cal{A}}}
\def\tUAres{\wt{U}_{\cal{A}}^{\mrm{res}}}
\def\tUfin{\wt{U}_{\vep}^{\mrm{fin}}}
\def\tUres{\wt{U}_{\vep}^{\mrm{res}}}
\def\tVdia{\wt{V}^{\mrm{dia}}}
\def\tVAres{\wt{V}_{\AA}^{\mrm{res}}}
\def\tVfin{\wt{V}^{\mrm{fin}}_{\vep}}
\def\tVnil{\wt{V}^{\mrm{nil}}_{\vep}}
\def\tVres{\wt{V}^{\mrm{res}}_{\vep}}
\def\UA{U_{\cal{A}}}
\def\UAres{U_{\cal{A}}^{\mrm{res}}}
\def\Ufin{U_{\varepsilon} ^{\mrm{fin}}}
\def\Ures{U_{\varepsilon}^{\mrm{res}}}
\def\VAres{V_{\AA}^{\mrm{res}}}
\def\Vfin{V^{\mrm{fin}}_{\vep}}
\def\Vnil{V^{\mrm{nil}}_{\vep}}
\def\Vres{V^{\mrm{res}}_{\vep}}
\def\tVq{\wt{V}_q}
\def\tPq{\wt{P}_q}
\def\tUq{\wt{U}_q}
\def\Pq{P_q}
\def\Vq{V_q}
\def\Uq{U_q}
\section{Introduction}
Let $\frg$ be a finite-dimensional complex simple Lie algebra 
and let $\wt{\frg}$ be the loop algebra of $\frg$. 
Let $q$ be an indeterminate 
and let $\Uq(\wt{\frg})$ be the quantum algebra over $\bbC(q)$ 
associated with $\wt{\frg}$ and $q$, 
where $\bbC(q)$ is the rational function field. 
Let $n$ be the rank of $\frg$ and let $I:=\{1, 2, \cdots, n\}$ 
be the index set.
For $\mbf{a} \in \bbC(q)^{\times}$ and an index $\xi \in I$, 
there exists a finite-dimensional irreducible 
representation of $\Uq(\wt{\frg})$ 
which is called a \it{fundamental representation} 
denoted by $\tVq(\pi_{\xi}^{\mbf{a}})$ (see \S 4.4). 

In 1997, Akasaka and Kashiwara gave the condition for 
$\tVq(\pi_{\xi_1}^{\mbf{a}_1}) \otimes \cdots 
\otimes \tVq(\pi_{\xi_r}^{\mbf{a}_r})$ 
to be irreducible in the case that 
$\wt{\frg}$ is of type $A_n^{(1)}$ and $C_n^{(1)}$ (see \cite{AK}).
In 2002, 
Varagnolo and Vasserot showed that condition 
in the simply laced case (see \cite{VV}) 
and Kashiwara gave that condition for arbitrary $\wt{\frg}$ 
(see \cite{K}). 
Moreover, in 2002, Chari showed
the sufficient condition for the tensor product 
of the finite-dimensional irreducible representations 
of $\Uq(\wt{\frg})$ to be a highest-weight representation 
(see \cite{C}). 

In particular, if $\wt{\frg}=A_n^{(1)}$,
we explicitly obtain the necessary and sufficient conditions 
for $\tVq(\pi_{\xi_1}^{\mbf{a}_1}) \otimes \cdots 
\otimes \tVq(\pi_{\xi_r}^{\mbf{a}_r})$ to be irreducible.
Indeed, we have the following theorem 
(see \cite{AK} and Theorem \ref{thm Main theorem-GQLA} of this paper).
\newtheorem*{thm main1}{Theorem}
\begin{thm main1}
Let $m \in \bbN$, $\xi_1, \cdots, \xi_m \in I$, 
and $\mbf{a}_1, \cdots, \mbf{a}_m \in \bbC(q)^{\times}$. 
The following conditions (a) and (b) are equivalent. 

(a) 
$\tVq(\pi_{\xi_1}^{\mbf{a}_1}) \otimes 
\cdots \otimes \tVq(\pi_{\xi_m}^{\mbf{a}_m})$ 
is an irreducible representation of $\Uq(A_n^{(1)})$.

(b) 
For any $1 \leq k \neq k^{'} \leq m$ 
and $1 \leq t \leq 
\mrm{min}(\xi_k, \xi_{k^{'}}, n+1-\xi_k, n+1-\xi_{k^{'}})$, 
\begin{eqnarray*}
  \frac{\mbf{a}_{k^{'}}}{\mbf{a}_k}
  \neq q^{\pm (2t+|\xi_k-\xi_{k^{'}}|)}. 
\end{eqnarray*}
\end{thm main1}

We want to extend this theorem 
for the restricted quantm algebras of type $A_n^{(1)}$ 
at roots of unity. 

Let $l$ be an odd integer greater than $3$, 
let $\vep$  be a primitive $l$-th root of unity, 
and let $\Ures(\wt{\frg})$ be the restricted quantum algebra 
over $\bbC$ associated with $\wt{\frg}$ and $\vep$ 
(see \cite{L89}, \cite{CP97}, and \S 6.1). 
For a nonzero complex number $\mbf{a}$ and an index $\xi \in I$, 
there exists a finite-dimensional irreducible 
representation of $\Ures(\wt{\frg})$ 
which is called a \it{fundamental representation} 
denoted by $\tVres(\pi_{\xi}^{\mbf{a}})$ (see \S 6.5). 

In 1997, Chari and Pressley showed that 
for any finite-dimensional irreducible 
$\Ures(\wt{\frg})$-representation $V$, 
there exist some nonzero complex numbers 
$\mbf{a}_1$, $\cdots$, $\mbf{a}_r$ 
and indexes $\xi_1$, $\cdots$, $\xi_r \in I$ 
such that $V$ is isomorphic to a subquotient of 
$\tVres(\pi_{\xi_1}^{\mbf{a}_1}) \otimes \cdots 
\otimes \tVres(\pi_{\xi_r}^{\mbf{a}_r})$.
However, the conditions for the irreducibility have not been given yet. 

So we consider the conditions in the case that 
$\wt{\frg}$ is of type $A_n^{(1)}$. 
The main theorem is as follows 
(see Theorem \ref{thm Main theorem-RQLA}): 
\newtheorem*{thm main2}{Theorem}
\begin{thm main2} 
Let $m \in \bbN$, $\xi_1, \cdots, \xi_m \in I$, 
and $\mbf{a}_1, \cdots, \mbf{a}_m \in \bbC^{\times}$. 
The following conditions (a) and (b) are equivalent. 

(a) 
$\tVres(\pi_{\xi_1}^{\mbf{a}_1}) \otimes 
\cdots \otimes \tVres(\pi_{\xi_m}^{\mbf{a}_m})$ 
is an irreducible representation of $\Ures(A_n^{(1)})$.

(b) 
For any $1 \leq k \neq k^{'} \leq m$ 
and $1 \leq t \leq 
\mrm{min}(\xi_k, \xi_{k^{'}}, n+1-\xi_k, n+1-\xi_{k^{'}})$, 
\begin{eqnarray*}
  \frac{\mbf{a}_{k^{'}}}{\mbf{a}_k}
  \neq \vep^{\pm (2t+|\xi_k-\xi_{k^{'}}|)}. 
\end{eqnarray*}
\end{thm main2}

The organization of this paper is as follows.
In \S 2, we fix some notations. 
In \S 3, we review the generic quantum algebras 
of type $A_n$ and $A_n^{(1)}$. 
In \S 4, we introduce the fundamental representations 
of the generic quantum algebras of type $A_n$ and $A_n^{(1)}$.  
In \S 5 (resp. \S 6, \S 7), we prove the main theorem 
for the generic quantum algebras of type $A_n^{(1)}$ 
(resp. the restricted quantum algebras of type $A_n^{(1)}$ 
at roots of unity, the small quantum algebras of type $A_n^{(1)}$). 
\section{Notations}
\setcounter{equation}{0}
\renewcommand{\theequation}{\thesection.\arabic{equation}}
We fix the following notations 
(see \cite{Kac}, \cite{BK}). 
Let $\sl_{n+1}$ be the finite-dimensional simple Lie algebra 
over $\bbC$ of type $A_n$. 
We define $I:=\{1,2, \cdots , n\}$.
Let $(\fra_{i,j})_{i,j \in I}$ 
be the Cartan matrix of $\sl_{n+1}$, 
that is, $\fra_{i,i}=2$, $\fra_{i,j}=-1$ if $|i-j|=1$, 
and $\fra_{i,j}=0$ otherwise. 
Let $\{\al_i\}_{i \in I}$ (resp. $\{\al_i^{\vee}\}_{i \in I}$) 
be the set of the simple roots (resp. simple coroots) 
of $\sl_{n+1}$ 
and let $\De$ (resp. $\De_{+}$) be the root system 
(resp. the set of positive roots) of $\sl_{n+1}$. 
Let $\frh=\bigoplus_{i \in I} \bbC \al_i^{\vee}$ 
be the Cartan subalgebra of $\sl_{n+1}$ 
and let $\frh^{*}=\bigoplus_{i \in I} \bbC \al_i$ 
be the $\bbC$-dual space of $\frh$. 
We have a $\bbC$-bilinear map 
$\langle, \rangle: \frh^{*} \times \frh \arr \bbC$ 
such that $\langle \al_j, \al_i^{\vee} \rangle =\fra_{i,j}$ 
for any $i, j \in I$. 
Let $Q:=\bigoplus_{i \in I}\bbZ \al_i$ 
(resp. $Q_{+}:=\bigoplus_{i \in I} \bbZ_{+} \al_i$) 
be the root lattice (resp. positive root lattice) of $\sl_{n+1}$, 
where $\bbZ_{+}:=\{0, 1, 2, \cdots\}$. 
Let $\{\La_i\}_{i \in I}$ be the fundamental weights of $\sl_{n+1}$, 
that is, 
\begin{eqnarray*}
\La_i:=\frac{1}{n+1}\{(n-i+1)\sum_{k=1}^i k\al_k
+i\sum_{k=i+1}^n(n-k+1)\al_k\} \in \frh^{*},
\label{def fundamental weights}
\end{eqnarray*}
(see [H], \S 13). 
We have 
$\langle \La_i, \al_j^{\vee} \rangle=\de_{i,j}$ for any $i,j \in I$. 
Let $P:=\bigoplus_{i \in I}\bbZ \La_i$ 
(resp. $P_{+}:=\bigoplus_{i \in I} \bbZ_{+} \La_i$) 
be the weight lattice (resp. positive weight lattice) 
of $\sl_{n+1}$.
Define a partial order $<$ in $P$ whereby 
\begin{eqnarray}
  \rm{$\nu \leq \nu^{'}$ if and only if $\nu^{'}-\nu \in Q_{+}$ 
  \q for $\nu, \nu^{'} \in P$}. 
\label{def root order} 
\end{eqnarray}
Let $\wt{\sl}_{n+1}
=\sl_{n+1} \otimes \bbC[t,t^{-1}]$ 
be the loop algebra of $\sl_{n+1}$. 
We define $\wt{I}:=I \sqcup \{0\}$ and 
\begin{eqnarray*}
 \fra_{0,0}:=2, \q \fra_{i,0}:=\fra_{0,j}:=0, 
\q \fra_{n,0}:=\fra_{0,n}:=-1 
\q \rm{for} \q 1 \leq i,j < n.
\end{eqnarray*}
Then $(\fra_{i,j})_{i,j \in \wt{I}}$ 
is the generalized Cartan matrix of $\wt{\sl}_{n+1}$. 
Let $\{\al_i\}_{i \in \wt{I}}$ be
the set of the simple roots of $\wt{\sl}_{n+1}$. 
We define $\wt{\frh}^{*}:=\bbC \al_0 \oplus \frh^{*}$.
We have a symmetric $\bbC$-bilinear form 
$(,): \wt{\frh}^{*} \times \wt{\frh}^{*} \arr \bbC$ such that 
$(\al_i, \al_j)= \fra_{i,j}$ for any $i,j \in \wt{I}$. 
Let $s_{i}$ be the simple reflection on $\wt{\frh}^{*}$, that is, 
$s_{i}(\la)=\la- (\la, \al_i)\al_i$ for $\la \in \wt{\frh}^{*}$. 
The affine Weyl group $\wt{\WW}$ of $\wt{\sl}_{n+1}$ 
(resp. Weyl group $\WW$ of $\sl_{n+1}$) 
is generated by $\{s_i\}_{i \in \wt{I}}$ (resp. $\{s_i\}_{i \in I}$). 

Let $q$ be an indeterminate. 
For $r \in \bbZ$ and $m \in \bbN:=\{1, 2, \cdots\}$, 
we define $q$-integers and Gaussian binomial coefficients 
in the rational function field $\bbC(q)$ whereby 
\begin{eqnarray*}
&&[r]_q:=\frac{q^{r}-q^{-r}}{q-q^{-1}}, 
\q [m]_q!:=[m]_q[m-1]_q \cdots [1]_q, 
\q \left[
\begin{array}{c}
r\\
m
\end{array}
\right]_q
:=\frac{[r]_q[r-1]_q \cdots [r-m+1]_q}{[1]_q[2]_q \cdots [m]_q}.
\end{eqnarray*}
Similarly, for $c \in \bbC \, (c \neq 0, \pm 1)$,  we define 
\begin{eqnarray*}
&&[r]_c:=\frac{c^{r}-c^{-r}}{c-c^{-1}}, 
\q [m]_c!:=[m]_c[m-1]_c \cdots [1]_c, 
\q \left[
\begin{array}{c}
r\\
m
\end{array}
\right]_c
:=\frac{[r]_c[r-1]_c \cdots [r-m+1]_c}{[1]_c[2]_c \cdots [m]_c}.
\end{eqnarray*}
We define $[0]_q!:=[0]_c!:=1$.
\section{Quantum algebras : the generic case}
\setcounter{equation}{0}
\renewcommand{\theequation}{\thesection.\arabic{equation}}
\subsection{Definitions and properties}
\begin{df}
\label{def QA}
Let $\tUq:=U_q(\wt{\sl}_{n+1})$ (resp. $U_q:=U_q(\sl_{n+1})$) 
be the associative $\bbC(q)$-algebra generated by 
$\{E_{i}, F_{i}, K_i^{\pm 1} \, | 
\, i \in \wt{I} \rm{ (resp. $i \in I$)} \}$  
with the following defining relations. 
We call $\tUq$ (resp. $\Uq$) the quantum algebra of type $A_n^{(1)}$ 
(resp. $A_n$) or quantum loop algebra of type $A_n$:    
\begin{eqnarray*}
&& K_i K_i^{-1}= K_i^{-1} K_i=1, \q K_i K_j=K_j K_i, 
\q K_0=\prod_{i \in I}K_i^{-1}, \\
&& K_i E_{j} K_i^{-1}=q^{\fra_{i,j}}E_{j},  
\q K_i F_{j} K_i^{-1}=q^{-\fra_{i,j}}F_{j}, \\ 
&& E_{i} F_{j}- F_{j} E_{i} = \de_{i,j} 
\frac{K_i-K_i^{-1}}{q-q^{-1}}, \\
&& \sum_{p=0}^{1-\fra_{i,j}} (-1)^p E_{i}^{(p)} E_{j}
E_{i}^{(1-\fra_{i,j}-p)}=
\sum_{p=0}^{1-\fra_{i,j}} (-1)^p F_{i}^{(p)} F_{j} 
F_{i}^{(1-\fra_{i,j}-p)}=0
 \q i \neq j, 
\end{eqnarray*}
for $i, j \in \wt{I}$ (resp. $i, j \in I$), where 
\begin{eqnarray}
E_{i}^{(m)}:= \displaystyle \frac{1}{[m]_q!} E_{i}^m, 
\q F_{i}^{(m)} := \displaystyle \frac{1}{[m]_q!} F_{i}^m 
\q \rm{for} \q m \in \bbN. 
\label{def divided power}
\end{eqnarray}
\end{df}

Let $U_q^{+}$  (resp. $U_q^{-}$, $U_q^{0}$ ) 
be the $\bbC(q)$-subalgebra of $U_q$ 
generated by $\{E_{i}\}_{i \in I}$ 
(resp. $\{F_{i}\}_{i \in I}$, $\{K_i^{\pm 1}\}_{i \in I}$). 
Then $\Uq$ has the triangular decomposition, that is, 
the multiplication defines an isomorphism of $\bbC(q)$-vector spaces: 
\begin{eqnarray}
  \Uq^{-} \otimes \Uq^0 \otimes \Uq^{+} \wt{\arr} \Uq,
\label{fac TD-GQA}
\end{eqnarray}
(see \cite{L90a}). 
It is well known that 
$\tUq$ (resp. $U_q$) has a Hopf algebra structure. 
The comultiplication $\De_H$, counit $\ep_H$, 
and antipode $S_H$ of $\tUq$ (resp. $U_q$) are given by 
\begin{eqnarray}
&& \De_H(E_i)=E_i \otimes 1 + K_i \otimes E_i, 
   \q \De_H(F_i)=F_i \otimes K_i^{-1} +1 \otimes F_i, 
   \q \De_H(K_i)=K_i \otimes K_i, 
   \label{fac DeH}\\
&& \ep_{H}(E_i)=\ep_H(F_i)=0, \q \ep_H(K_i)=1, 
   \label{fac epH}\\
&& S_H(E_i)=-K_i^{-1}E_i, \q S_H(F_i)=-F_iK_i, 
   \q S_H(K_i)=K_i^{-1}, 
\label{fac SH}
\end{eqnarray} 
where $i \in \wt{I}$ (resp. $I$) 
(see \cite{Ja}, \S 3--\S 4 and \cite{CP95}, \S 2). 

Let $T_i$ be the $\bbC(q)$-algebra automorphism of $\tUq$ 
(resp. $\Uq$) introduced by Lusztig in \cite{L93}, \S 37: 
\begin{eqnarray*}
&& T_i(E_i^{(m)})=(-1)^m q^{-m(m-1)}F_i^{(m)}K_i^m, 
\q T_i(F_i^{(m)})=(-1)^m q^{m(m-1)}K_i^{-1}E_i^{(m)},  \no\\
&& T_i(E_j^{(m)})=\sum_{p=0}^{-m \fra_{i,j}}(-1)^{p-m\fra_{i,j}}
q^{-p}E_i^{(-m\fra_{i,j}-p)}E_j^{(m)}E_i^{(p)} \q (i \neq j), \no \\
&& T_i(F_j^{(m)})=\sum_{p=0}^{-m \fra_{i,j}}(-1)^{p-m\fra_{i,j}}
q^{p}F_i^{(p)}F_j^{(m)}F_i^{(-m\fra_{i,j}-p)} \q (i \neq j), \no\\
&& T_i(K_j)=K_jK_i^{-\fra_{i,j}}, 
\label{pro RA}
\end{eqnarray*}
where $i, j \in \wt{I}$ (resp. $i, j \in I$) and $m \in \bbN$. 
For $w \in \wt{\WW}$, 
let $w=s_{i_1} \cdots s_{i_m} 
(i_1, \cdots , i_m \in \wt{I}, m \in \bbN)$ 
be a reduced expression of $w$. 
Then $T_{w}:=T_{i_1} \cdots T_{i_m}$ 
is a well-defined automorphism of $\tUq$, that is, 
$T_{w}$ does not depend on the choice of the reduced expression of $w$.

It is known that $\tUq$ has another realization 
which is called the Drinfel'd realization. 
\begin{df}[\cite{D}]
\label{def DR} 
Let $\DD_q$ be an associative $\bbC(q)$-algebra 
generated by 
$\{X_{i,r}^{\pm}, H_{i,s}, K_i^{\pm 1} \, | 
\, i \in I, r, s \in \bbZ, s \neq 0\}$ 
with the defining relations 
\begin{eqnarray*}
&& K_iK_i^{-1}=K_i^{-1}K_i=1, 
\q [K_i, K_j]=[K_i, H_{j,s}]=[H_{i,s}, H_{j,s^{'}}]=0, \\
&& K_iX_{j,r}^{\pm}K_i^{-1}=q^{\pm \fra_{i,j}}X_{j,r}^{\pm}, 
\qq [H_{i,s}, X_{j,r}^{\pm}]=\pm \frac{[s\fra_{i,j}]_q}{s}X_{j,r+s}^{\pm}, \\
&&X_{i,r+1}^{\pm}X_{j,r^{'}}^{\pm}
-q^{\pm \fra_{i,j}}X_{j,r^{'}}^{\pm}X_{i,r+1}^{\pm}
=q^{\pm \fra_{i,j}}X_{i,r}^{\pm}X_{j,r^{'}+1}^{\pm}
-X_{j,r^{'}+1}^{\pm}X_{i,r}^{\pm}, \\
&& [X_{i,r}^{+}, X_{j,r^{'}}^{-}]
=\de_{i,j}\frac{\Psi_{i,r+r^{'}}^{+}-\Psi_{i,r+r^{'}}^{-}}{q-q^{-1}}, \\
&&\sum_{\pi \in \SS_{\bf{m}}} \sum_{p=0}^{\bf{m}} (-1)^p 
\left[
\begin{array}{r}
\bf{m}\\
p
\end{array}
\right]_q
X_{i,r_{\pi(1)}}^{\pm} \cdots X_{i,r_{\pi(p)}}^{\pm} 
X_{j,r^{'}}^{\pm} X_{i,r_{\pi(p+1)}}^{\pm} 
\cdots X_{i,r_{\pi(\bf{m})}}^{\pm}=0, 
\q (i \neq j),
\end{eqnarray*}
for $r_1, \cdots, r_{\bf{m}} \in \bbZ$, 
where $\bf{m}:=1-\fra_{i,j}$, 
$\SS_{\bf{m}}$ is the symmetric group on $\bf{m}$ letters, 
and $\Psi_{i,r}^{\pm}$ are determined by 
\begin{eqnarray*}
\sum_{r=0}^{\infty}\Psi_{i,\pm r}^{\pm}u^{\pm r}
:=K_i^{\pm 1}\mrm{exp}
(\pm (q-q^{-1})\sum_{s=1}^{\infty}H_{i,\pm s}u^{\pm s}) 
\q \rm{in} \q \DD_q[[u]],
\end{eqnarray*}
and $\Psi_{i,\pm r}^{\pm}:=0$ if $r<0$.
\end{df}
\begin{thm}[\cite{B}]
\label{thm DR} 
  $\DD_q$ is isomorphic to $\tUq$ as a $\bbC(q)$-algebra.
\end{thm}
The isomorphism from $\DD_q$ to $\tUq$ is given in \cite{B}. 
Here we introduce an isomorphism from $\tUq$ to $\DD_q$ 
introduced in \cite{CP94a}, \S 2.5 (see also \cite{AN}, II-C).
There exists a $\bbC(q)$-algebra isomorphism 
$T : \tUq \arr \DD_q$ such that 
\begin{eqnarray*}
&& T(E_i)=X_{i,0}^{+}, \q T(F_i)=X_{i,0}^{-}, 
   \q T(K_i)=K_i, \no \\
&& T(E_0)=(-1)^{m+1}q^{n+1}[X_{n,0}^{-},  \cdots [
   X_{m+1,0}^{-},[ X_{1,0}^{-},  \cdots [X_{m-1,0}^{-},
   X_{m,1}^{-}]_{q^{-1}} \cdots ]_{q^{-1}}\prod_{i \in I}K_i^{-1}, \no \\
&& T(F_0)=(-1)^{m+n}[X_{n,0}^{+},  \cdots [
   X_{m+1,0}^{+},[ X_{1,0}^{+},  \cdots [X_{m-1,0}^{+},
   X_{m,-1}^{+}]_{q^{-1}} \cdots ]_{q^{-1}}\prod_{i \in I} K_i, 
\label{pro DR}
\end{eqnarray*}
for $m, i \in I$, 
where $[u,v]_{q^{\pm1}}:=uv-q^{\pm1}vu$ for $u,v \in \tUq$ 
($T$ is independent of the choice of $m$). 
We identify $\tUq$ with $\DD_q$ by this isomorphism.

 Now, for $i \in I$ and $r \in \bbZ$, we define $\PP_{i,r} \in \tUq$ whereby 
\begin{eqnarray}
  \sum_{m=1}^{\infty} \PP_{i, \pm m}u^m:=
  \mrm{exp}(-\sum_{s=1}^\infty 
  \frac{q^s}{[s]_q}H_{i,s}u^s), 
  \q \PP_{i,0}:=1,
  \q \rm{in} \q \tUq[[u]].
\label{def PP}
\end{eqnarray}

Let $\tUq^{\pm}$  (resp. $\tUq^{0}$) 
be the $\bbC(q)$-subalgebra of $\tUq$ generated by 
$\{X_{i,r}^{\pm} \, | \, i \in I, r \in \bbZ\}$ 
(resp. $\{K_i^{\pm 1}, \PP_{i,r} \, | \, i \in I, r \in \bbZ\}$). 
Then $\tUq$ also has the triangular decomposition, that is, 
the multiplication defines an isomorphism of $\bbC(q)$-vector spaces: 
\begin{eqnarray}
  \tUq^{-} \otimes \tUq^0 \otimes \tUq^{+} \wt{\arr} \tUq,
\label{fac TD-GQLA}
\end{eqnarray}
(see \cite{CP01}, \S 4, \cite{BCP}, and \cite{B}).
We define 
\begin{eqnarray*}
 X_{\pm}:=\sum_{i \in I, r \in \bbZ} \bbC(q)X_{i,r}^{\pm}, 
\q X_{\pm}(i):=\sum_{j \in (I\bs\{i\}), r \in \bbZ}
\bbC(q)X_{j,r}^{\pm}, 
\label{def X+-}
\end{eqnarray*}
where $X_{\pm}(i):=0$ if $I=\{i\}$.
\begin{pro}[\cite{C}, Proposition 2.2]
\label{pro C02-2.2}
Let $i \in I$, $r \in \bbZ$, and $s \in \bbZ^{\times}$. 

(a) Modulo $\tUq X_{-} \otimes \tUq (X_{+})^2
+ \tUq X_{-} \otimes \tUq X_{+}(i)$,
\begin{eqnarray*}
&& \De_H(X_{i,r}^{+})
   =X_{i,r}^{+} \otimes 1 + K_i \otimes X_{i,r}^{+}
   + \sum_{p=1}^{r}\Psi_{i,p}^{+} \otimes X_{i,r-p}^{+} 
   \q \rm{if} \q r \geq 0, \\
&& \De_H(X_{i,-r}^{+})
   = K_i^{-1} \otimes X_{i,-r}^{+}+X_{i,-r}^{+} \otimes 1 
   + \sum_{p=1}^{r-1}\Psi_{i,-p}^{-} \otimes X_{i,-r+p}^{+} 
   \q \rm{if} \q r > 0. 
\end{eqnarray*}

(b) Modulo $\tUq (X_{-})^2 \otimes \tUq X_{+}
+ \tUq X_{-} \otimes \tUq X_{+}(i)$,
\begin{eqnarray*}
&& \De_H(X_{i,r}^{-})
   =X_{i,r}^{-} \otimes K_i + 1 \otimes X_{i,r}^{-}
   + \sum_{p=1}^{r-1} X_{i,r-p}^{-}  \otimes \Psi_{i,p}^{+}
   \q \rm{if} \q r > 0, \\
&& \De_H(X_{i,-r}^{-})
   = X_{i,-r}^{+} \otimes K_i^{-1} + 1 \otimes X_{i,-r}^{-}  
   + \sum_{p=1}^{r}X_{i,-r+p}^{-} \otimes \Psi_{i,-p}^{-} 
   \q \rm{if} \q r \geq 0. 
\end{eqnarray*}

(c) Modulo $\tUq X_{-} \otimes \tUq X_{+}$,
\begin{eqnarray*}
  \De_H(H_{i,s})=H_{i,s} \otimes 1 + 1 \otimes H_{i,s}.
\end{eqnarray*}
\end{pro}
\subsection{Representation theory of $\Uq$ and $\tUq$}
\begin{df} 
\label{def HWM}
Let $V$ be a representation of $\tUq$ (resp. $\Uq$). 

(i) Let $v \in V$. 
If $X_{i,r}^{+}v=0$ for all $i \in I$, $r \in \bbZ$ 
(resp. $E_iv=0$ for all $i \in I$), 
we call $v$ a \it{pseudo-primitive vector} 
(resp. \it{primitive vector}) in $V$. 

(ii) For any $\mu \in P$, we define 
\begin{eqnarray*}
V_{\mu}:=\{v \in V \, | \, K_iv=q^{\langle \mu, \al_i^{\vee} \rangle}v 
\q \rm{for all $i \in I$}\}. 
\end{eqnarray*}
If $V_{\mu} \neq 0$, we call $V_{\mu}$ a \it{weight space} of $V$. 
For $v \in V_{\mu}$, we call $v$ a \it{weight vector} with weight $\mu$ 
and define $\mrm{wt}(v):=\mu$.

(iii) For any $\bbC$-algebra homomorphism $\La: \tUq^0 \arr \bbC(q)$, 
we define 
\begin{eqnarray*}
V_{\La}:=\{v \in V \, | \, uv=\La(u)v
\q \rm{for all $u \in \tUq^0$}\}. 
\end{eqnarray*}
If $V_{\La} \neq 0$, we call $V_{\La}$ a \it{pseudo-weight space} of $V$. 
For $v \in V_{\La}$, we call $v$ 
a \it{pseudo-weight vector} with \it{pseudo-weight} $\La$ 
and define $\mrm{pwt}(v):=\La$.

(iv) Let $\La: \tUq^0 \arr \bbC$ be a $\bbC$-algebra homomorphism 
and $\la$ be an element in $P_{+}$. 
If there exists a nonzero pseudo-primitive vector $v_{\La} \in V_{\La}$ 
(resp. primitive vector $v_{\la} \in V_{\la}$)
such that $V=\tUq v_{\La}$ (resp. $V=\Uq v_{\la}$), 
we call $V$ a \it{pseudo-highest weight representation} of $\tUq$
(resp. \it{highest-weight representation} of $\Uq$) 
with the \it{pseudo-highest weight $\La$} 
(resp. \it{highest weight $\la$})
generated by a \it{pseudo-highest weight vector $v_{\La}$} 
(resp. \it{highest-weight vector $v_{\La}$}). 
\end{df}

Let $V$ be a representation of $\tUq$ (resp. $\Uq$). 
We call $V$ of \it{type $\mbf{1}$} if
\begin{eqnarray*}
V=\bigoplus_{\mu \in P} V_{\mu}.
\label{def type 1}
\end{eqnarray*}
In general, finite-dimensional representations of $\tUq$ (resp. $\Uq$) 
are classified into $2^n$ types according to 
$\{\sigma : Q \arr \{\pm 1\}; \rm{ group homomorphism}\}$. 
It is known that for any $\sigma : Q \arr \{\pm 1\}$, 
the category of finite-dimensional representations 
of $\tUq$ (resp. $\Uq$) of type $\sigma$ is essentially equivalent to 
the category of the finite-dimensional representations 
of $\tUq$ (resp. $\Uq$) of type $\mbf{1}$ 
(see \cite{CP94b}, \S 10--\S 11). 

For any $\tUq$-representation (resp. $\Uq$-representation) V, 
we have
\begin{eqnarray}
  X_{i,r}^{\pm} V_{\mu} \subset V_{\mu \pm \al_i}, 
  \q (\rm{resp.} \q E_i V_{\mu} \subset V_{\mu+\al_i}, 
  \q F_i V_{\mu} \subset V_{\mu-\al_i}), 
\label{fac weight space}
\end{eqnarray}
where $i \in I$, $r \in \bbZ$, and $\mu \in P$. 

The following theorem is well known (see \cite{Ja}).
\begin{thm} 
\label{thm CT-GQA}
For any $\la \in P_{+}$, 
there exists a unique (up to isomorphism) 
finite-dimensional irreducible highest-weight $\Uq$-representation 
$\Vq(\la)$ with the highest weight $\la$ of type $\mbf{1}$. 
Conversely, for any finite-dimensional irreducible 
$\Uq$-representation $V$ of type $\mbf{1}$, 
there exists a unique $\la \in P_{+}$
such that $V$ is isomorphic to $\Vq(\la)$ 
as a representation of $\Uq$.
\end{thm}
We define a set of polynomials $\bbC(q)_0[t]$ whereby
\begin{eqnarray*}
\bbC(q)_0[t]:=\{\pi(t) \in \bbC(q)[t] \, | \, 
\rm{there exist some $a_1, \cdots, a_r \in \bbC(q)$ 
such that $\pi(t)=\prod_{s=1}^r(1-a_st)$}\}. 
\label{def bbC0}
\end{eqnarray*}
For any $\mbf{\pi}=(\pi_i(t))_{i \in I} \in (\bbC(q)_0[t])^n$, 
there exists a unique $\bbC(q)$-algebra homomorphism 
$\La_{\mbf{\pi}}: \tUq^0 \arr \bbC(q)$ such that 
\begin{eqnarray*}
  \La_{\mbf{\pi}}(K_i^{\pm 1})=q^{\pm \mrm{deg}\pi_i(t)}, \q
  \sum_{m=1}^{\infty} \La_{\mbf{\pi}}(\PP_{i,\pm m})t^m 
  =\pi_i^{\pm}(t), 
  \label{def Lapi}
\end{eqnarray*}
where 
\begin{eqnarray}
\pi_i^{+}(t):=\pi_i(t), 
\q \pi_i^{-}(t):=\frac{t^{\mrm{deg}\pi_i(t)}\pi_i(t^{-1})}
{(t^{\mrm{deg}\pi_i(t)}\pi_i(t^{-1}))|_{t=0}}.
\label{def pi-}
\end{eqnarray}
From \cite{CP97}, \S 3, for $i \in I$, we have
\begin{eqnarray}
  \sum_{m=0}^{\infty} \La_{\mbf{\pi}}(\Psi_{i,m}^{+}) u^m=
  q^{\mrm{deg}(\pi_i(u))}\frac{\pi_i(q^{-2}u)}{\pi_i(u)}
  =\sum_{m=0}^{\infty} \La_{\mbf{\pi}}(\Psi_{i,-m}^{-}) u^{-m} 
  \q \in \tUq[[u]],
\label{fac Psi}
\end{eqnarray}
in the sense that left- and right-hand terms are 
the Laurent expansions of the middle term about $0$ and $\infty$, 
respectively 
(see also \cite{CP95}, Theorem 3.3 and \cite{CP94b}, 12.2B). 

For any pseudo-highest weight representation of $\tUq$ 
with the pseudo-highest weight $\La_{\mbf{\pi}}$, 
we simply call it a pseudo-highest representation of $\tUq$ 
with the pseudo-highest weight $\mbf{\pi}$. 
\begin{thm}[\cite{CP97}, \S 2 and \cite{C}, \S 2]
\label{thm CT-GQLA}
For any $\mbf{\pi} \in (\bbC(q)_0[t])^n$, 
there exists a unique (up to isomorphism) 
finite-dimensional irreducible pseudo-highest weight 
$\tUq$-representation $\tVq(\mbf{\pi})$ 
with the pseudo-highest weight $\mbf{\pi}$ of type $\mbf{1}$. 
\end{thm}
From (\ref{fac TD-GQA}), (\ref{fac TD-GQLA}), 
and (\ref{fac weight space}), 
for $\la \in P_{+}$ and 
$\mbf{\pi}=(\pi_i(t))_{i \in I} \in (\bbC(q)_0[t])^n$, we have
\begin{eqnarray}
 \mrm{dim}_{\bbC(q)}(\Vq(\la)_{\la})=1, 
 \q \Vq(\la)=\bigoplus_{\nu \leq \la} \Vq(\la)_{\nu},
 \q \mrm{dim}_{\bbC(q)}(\tVq(\mbf{\pi})_{\la(\mbf{\pi})})=1, 
 \q \tVq(\mbf{\pi})
 =\bigoplus_{\nu \leq \la(\mbf{\pi})} \tVq(\mbf{\pi})_{\nu}, 
\label{fac weight decomposition}
\end{eqnarray}
where $<$ is the partial order defined in (\ref{def root order}) 
and 
\begin{eqnarray}
 \la(\mbf{\pi}):=\sum_{i \in I}\la_i(\pi)\La_i
 :=\sum_{i \in I}\mrm{deg}(\mbf{\pi}_i(t))\La_i \in P_{+}. 
\label{def lapi} 
\end{eqnarray}
It is known that for any $\om \in \WW$ and $\mu \in P$, 
\begin{eqnarray}
\label{fac weight}
 \mrm{dim}_{\bbC(q)}(\Vq(\la)_{\om \mu})
 =\mrm{dim}_{\bbC(q)}(\Vq(\la)_{\mu}).
\end{eqnarray}
\begin{pro}[\cite{CP95}, Proposition 3.4]
\label{pro tensor highest}
Let $V$ (resp. $V^{'}$) be a pseudo-highest weight representations 
of $\tUq$ with the pseudo-highest weight $\mbf{\pi}$ 
(resp. $\mbf{\pi}^{'}$) 
generated by a pseudo-highest weight vector 
$v_{\mbf{\pi}}$ (resp. $v_{\mbf{\pi}^{'}}^{'}$). 
Then $v_{\mbf{\pi}} \otimes v_{\mbf{\pi}^{'}}^{'}$ is 
a pseudo-primitive vector with 
$\mrm{pwt}(v_{\mbf{\pi}} \otimes v_{\mbf{\pi}^{'}}^{'})
=\La_{\mbf{\pi} \mbf{\pi}^{'}}$.
\end{pro}
\begin{cor}[\cite{CP95}, Corollary 3.5 and \cite{CP94b}, Theorem 12.2.6]
\label{cor tensor highest}
Let $\mbf{\pi}, \mbf{\pi}^{'} \in (\bbC(q)_0[t])^n$ and 
let $v_{\mbf{\pi}}$ (resp. $v_{\mbf{\pi}^{'}}^{'}$) 
be a pseudo-highest weight vector in 
$\tVq(\mbf{\pi})$ (resp. $\tVq(\mbf{\pi}^{'})$). 
$\tVq(\mbf{\pi}\mbf{\pi^{'}})$ is isomorphic to 
a quotient of the subrepresentation 
of $\tVq(\mbf{\pi}) \otimes \tVq(\mbf{\pi}^{'})$ 
generated by $v_{\mbf{\pi}} \otimes v_{\mbf{\pi}^{'}}^{'}$. 
In particular, 
if $\tVq(\mbf{\pi}) \otimes \tVq(\mbf{\pi}^{'})$ is irreducible, 
$\tVq(\mbf{\pi}) \otimes \tVq(\mbf{\pi}^{'})$ 
is isomorphic to 
$\tVq(\mbf{\pi}\mbf{\pi}^{'})$. 
\end{cor}
\subsection{The dual and involution representations of $\tUq$}
For any $\tUq$-representation $V$, 
let $V^{*}=\{\rm{$g: V \arr \bbC(q)$; $\bbC(q)$-linear map}\}$ 
be the dual $\tUq$-representation of $V$. 
The action of $\tUq$ on $V^{*}$ is defined by
\begin{eqnarray*}
 (u.g)(v):=g(S_H(u).v), 
 \q \rm{for $u \in \tUq$, $g \in V^{*}$, and $v \in V$},
\end{eqnarray*}  
where $S_H$ is as in (\ref{fac SH}).

There exists a $\bbC(q)$-algebra involution 
$\Om: \tUq \arr \tUq$ such that 
\begin{eqnarray}
 \Om(X_{i,r}^{\pm})=-X_{i,-r}^{\mp}, 
\q \Om(H_{i,s})=-H_{i,-s}, 
\q \Om(\Psi_{i,r}^{\pm})=\Psi_{i,-r}^{\mp}, 
\q \Om(K_i^{\pm 1})=K_i^{\mp 1},
\label{fac Om}
\end{eqnarray}
for $i \in I$, $r \in \bbZ$, and $s \in \bbZ^{\times}$ 
(see \cite{CP96}, Proposition 1.4(b)). 
For any $\tUq$-representation $V$, 
let $V^{\Om}$ be the pull-back of $V$ through $\Om$. 
For any finite-dimensional $\tUq$-representations 
$V$ and $W$, we obtain 
\begin{eqnarray}
 (V \otimes W)^{*} \cong W^{*} \otimes V^{*}, 
\q (V \otimes W)^{\Om} \cong W^{\Om} \otimes V^{\Om} 
\q \rm{as representations of $\tUq$}.
\label{pro DR+IR tensor}
\end{eqnarray} 
\begin{pro}
\label{pro DR+IR=irreducible}
Let $V$ be a finite-dimensional representation of $\tUq$. 
If $V^{*}$ and $V^{\Om}$ are pseudo-highest weight 
representations of $\tUq$, 
$V$ is irreducible as a representation of $\tUq$.
\end{pro}

Proof.
Let $v_H$ be a pseudo-highest weight vector in $V^{\Om}$ 
and let $\la$ be the weight of $v_H$. 
Since $V^{\Om}$ is a pseudo-highest weight 
representation,  
from (\ref{fac TD-GQLA}), we have 
\begin{eqnarray*}
 V^{\Om}=\tUq v_H=\bbC v_H \oplus(\bigoplus_{\mu<\la}V^{\Om}_{\mu}).
\end{eqnarray*}
Hence, from the definition of $V^{\Om}$, we have 
\begin{eqnarray*}
 V=\bbC v_H \oplus(\bigoplus_{\mu>-\la}V_{\mu}).
\end{eqnarray*}
Let $W$ be a proper $\tUq$-subrepresentation of $V$. 
We shall prove $W=0$.
Since $V$ is generated by $v_H$ as a representation of $\tUq$, 
$v_H$ is not included in $W$. 
Hence we have
\begin{eqnarray}
 W \subset (\bigoplus_{\mu>-\la}V_{\mu}). 
\label{pro DR+IR=irreducible 1} 
\end{eqnarray}

Now we define 
\begin{eqnarray*}
 \wt{W}:=\{g \in V^{*} \, | \, g|_W=0\},
\end{eqnarray*}
and let $g_H$ be an element in $V^{*}$ defined by 
\begin{eqnarray*}
 g_H(v_H):=1, \q g(\bigoplus_{\mu > -\la} V_{\mu}):=0.
\end{eqnarray*}
Then, from (\ref{pro DR+IR=irreducible 1}), 
$g_H$ is included in $\wt{W}$. 
Moreover the weight of $g_H$ is 
maximam in $U_q$-representation $V^{*}$. 
Then, since $V^{*}$ is a pseudo-highest weight $\tUq$-representation, 
$g_H$ is a pseudo-highest weight vector in $V^{*}$. 
Thus we have 
\begin{eqnarray*}
 V^{*}=\tUq g_H \subset \wt{W}.
\end{eqnarray*} 
Therefore $W=0$. 
\qed \\

\begin{pro}[\cite{CP96}, Proposition 5.1 and \cite{C}, (2.20), (2.21)] 
\label{pro DR+IR}
Let $\mbf{\pi}=(\pi_i(t))_{i \in I} \in (\bbC(q)_0[t])^n$. 

(a) There exists an integer $c \in \bbZ$ 
depending only on $\sl_{n+1}$ such that 
\begin{eqnarray*}
 V_q(\mbf{\pi})^{*} \cong 
V_q(\pi_n(q^{c} t), \cdots, \pi_1(q^{c} t)).
\end{eqnarray*} 

(b) There exists a nonzero complex number $\kappa \in \bbC^{\times}$ 
depending only on $\sl_{n+1}$ such that 
\begin{eqnarray*}
 V_q(\mbf{\pi})^{\Om} \cong 
V_q(\pi_n^{-}(q^2\kappa t), \cdots, \pi_1^{-}(q^2\kappa t)),
\end{eqnarray*} 
\end{pro}
where $\pi_i^{-}(t)$ is as in (\ref{def pi-}).
\subsection{The extremal vectors in $\tUq$}
In this subsection, 
we introduce the extremal vectors in $\tVq(\mbf{\pi})$ 
(see \cite{C}, \S 4 and \cite{AK}).
For $w \in W$, let $w=s_{i_1}s_{i_2} \cdots s_{i_k}$ 
be a reduced expression of $w$. 
For $\la=\sum_{i \in I}\la_i\La_i \in P_{+}$, 
we define 
\begin{eqnarray*}
 m_{i_j}^{\la}:= \langle s_{i_{j+1}}s_{i_{j+2}} 
 \cdots s_{i_k}\la, \al_{i_j}^{\vee} \rangle,
\end{eqnarray*}
where $m_{i_k}^{\la}:=\la_{i_k}$. 
Then we have $m_{i_j}^{\la} \geq 0$ for all $1 \leq j \leq k$. 
\begin{df} 
\label{def EV}
For $\mbf{\pi}=(\pi_i(t))_{i \in I} \in (\bbC(q)_0[t])^n$, 
let $v_{\mbf{\pi}}$ be a pseudo-highest weight vector 
in $\tVq(\mbf{\pi})$ 
and let $\la(\pi)$ be as in (\ref{def lapi}). 
For $w \in W$, let $w=s_{i_1}s_{i_2} \cdots s_{i_k}$ 
be a reduced expression of $w$. 
For $1 \leq j \leq k$, we define 
\begin{eqnarray*}
 v_{\mbf{\pi}}(s_{i_j}s_{i_{j+1}} \cdots s_{i_k})
:=F_{i_j}^{(m_{i_j}^{\la(\mbf{\pi})})}
F_{i_{j+1}}^{(m_{i_{j+1}}^{\la(\mbf{\pi})})} 
\cdots F_{i_k}^{(m_{i_k}^{\la(\mbf{\pi})})}v_{\mbf{\pi}},
\end{eqnarray*}
which is called a \it{extremal vector} in $\tVq(\mbf{\pi})$.
\end{df}
We have 
\begin{eqnarray}
  X_{i_j,r}^{+}
  v_{\mbf{\pi}}(s_{i_{j+1}}s_{i_{j+2}} \cdots s_{i_k})=0, 
\label{fac EV}
\end{eqnarray}
for any $1 \leq j \leq k$ and $r \in \bbZ$.
For $i, j \in I$ such that $i \leq j$, we define
\begin{eqnarray}
 \om_{i,j}:=
(s_i s_{i+1} \cdots s_j)(s_1s_2 \cdots s_{j-1})(s_1s_2 \cdots s_{j-2})
\cdots (s_1s_2)s_1. 
\label{def om}
\end{eqnarray}
Then $\om_{1,n}$ is the longest element in $\WW$. 
For any $\la=\sum_{i \in I} \la_i\La_i \in P_{+}$, we have 
\begin{eqnarray}
&& \langle \om_{i,j}\la, \al_{i-1}^{\vee} \rangle
   =\la_{j-i+2}+\la_{j-i+3}+ \cdots + \la_j, 
   \q \rm{if} \q i \neq 1, 
\label{fac swn1}\\
&& \langle \om_{1,j}\la, \al_{j+1}^{\vee} \rangle
   =\la_1+\la_2+ \cdots +\la_{j+1}.
\label{fac sw1} 
\end{eqnarray}
Hence we obtain 
\begin{eqnarray}
&& v_{\pi}(s_{i-1}\om_{i,j})
   =F_{i-1}^{(\la_{j-i+2}(\mbf{\pi})+ \cdots + \la_j(\mbf{\pi}))}
   v_{\pi}(\om_{i,j}), 
   \q \rm{if} \q i \neq 1,
\label{fac vomn1}\\
&& v_{\pi}(s_{j+1}\om_{1,j})
   =F_{j+1}^{(\la_1(\mbf{\pi})+ \cdots + \la_{j+1}(\mbf{\pi}))}
   v_{\pi}(\om_{1,j}). 
\label{fac vom1}
\end{eqnarray}
\begin{pro}[\cite{C}, Proposition 4.1 and Proposition 6.3 (i)] 
\label{pro EV}
For $\mbf{\pi} \in (\bbC(q)_0[t])^n$, 
let $v_{\mbf{\pi}}$ be a pseudo-highest vector in $\tVq(\mbf{\pi})$. 
For $i,j \in I$ such that $i \leq j$, 
we have 
\begin{eqnarray*}
&& H_{i-1,1}v_{\mbf{\pi}}(\om_{i,j})
   =\sum_{k=j-i+2}^j q^{2j-i-k+1}\La_{\mbf{\pi}}(H_{k,1})
   v_{\mbf{\pi}}(\om_{i,j}) 
   \q (i \neq 1), \\
&& H_{j+1,1}v_{\mbf{\pi}}(\om_{1,j})
   =\sum_{k=1}^{j+1} q^{j+1-k}\La_{\mbf{\pi}}(H_{k,1})
   v_{\mbf{\pi}}(\om_{1,j}).
\end{eqnarray*}
\end{pro}
\section{The fundamental representations: the generic case}
\setcounter{equation}{0}
\renewcommand{\theequation}{\thesection.\arabic{equation}}
\subsection{The fundamental representations of $\Uq$}
For $\la \in P_{+}$, 
we call $\la$ \it{minimal weight} 
if $\mu \in P_{+}$, $\mu \leq \la$ implies that $\mu=\la$. 
Moreover, we call a nonzero minimal weight 
a \it{minuscule weight} 
(see \cite{Ja}, \S 5A.1).
\begin{pro}[\cite{H}, \S 13, Exercises 13] 
\label{pro minimal}
Let $\la \in P_{+}$.

(a) $\la$ is minimal weight if and only if 
$\langle \la, \al^{\vee} \rangle \in \{0, \pm 1\}$ 
for all $\al \in \De$. 

(b) For $\xi \in I$, $\La_{\xi}$ is a minuscule weight. 
Conversely, if $\la$ is a minuscule weight, 
there exists an index $\xi \in I$ such that $\la=\La_{\xi}$. 
\end{pro}
For $\xi \in I$, we call $\Vq(\La_\xi)$ 
a \it{fundamental representation} of $\Uq$. 
We can construct these representations as follows.
Let $W_q(\La_{\xi})$ be a $\#(\WW \La_{\xi})$-dimensional 
$\bbC(q)$-vector space 
and let $\{z_{\mu} \, | \, \mu \in \WW \La_{\xi}\}$ 
be a $\bbC(q)$-basis of $W_q(\La_{\xi})$. 
\begin{pro}[\cite{Ja}, \S 5A.1]
\label{pro FR-GQA}
(a) We can define a $\Uq$-representation structure on $W_q(\La_{\xi})$ 
by the following formula: 
\begin{eqnarray}
&& K_i^{\pm 1}z_{\mu}
=q^{\pm (\mu, \al_i)}z_{\mu}, 
\label{pro FR-GCA-K} \\
&& E_iz_{\mu}=
\begin{cases}
z_{\mu+\al_i}, & \langle \mu, \al_i^{\vee} \rangle= -1, \\
0, & otherwise, 
\end{cases}
\label{pro FR-GCA-E} \\
&& F_iz_{\mu}=
\begin{cases}
z_{\mu-\al_i}, & \langle \mu, \al_i^{\vee} \rangle= 1, \\
0, & otherwise, 
\label{pro FR-GCA-F} 
\end{cases}
\end{eqnarray}
for any $i \in I$ and $\mu \in \WW\La_{\xi}$.

(b) $W_q(\La_{\xi})$ is isomorphic to $\Vq(\La_{\xi})$ 
as a representation of $\Uq$.
\end{pro}
We identify $W_q(\La_{\xi})$ with $V_q(\La_{\xi})$. 
Obviously, $z_{\La_{\xi}}$ is a highest-weight vector 
in $W_q(\La_{\xi})$. 
From (\ref{fac weight decomposition}) and (\ref{fac weight}), 
for any $\mu \in \WW \La_{\xi}$,
we have
\begin{eqnarray}
 \mrm{dim}_{\bbC(q)}(\Vq(\La_{\xi})_{\mu})=1. 
 \label{pro dim of FR}
\end{eqnarray}
\subsection{Another realization of 
the fundamental representations of $\Uq$}
We define $\wh{I}:=I \sqcup \{n+1\}=\{1 ,2, \cdots, n+1\}$.
For $\xi \in I$, we define 
\begin{eqnarray*}
 \JJ_\xi:=\{J=\{j_1, j_2, \cdots, j_\xi\} \subset \wh{I} \, 
 | \, j_1 < j_2 < \cdots < j_\xi\}.
\label{def JJ*}
\end{eqnarray*}
Let $L_q(\La_{\xi})$ be a $\#(\JJ_\xi)$-dimensional $\bbC(q)$-vector space. 
We regard $\JJ_\xi$ as a $\bbC(q)$-basis of $L_q(\La_{\xi})$. 
\begin{pro}[\cite{DO}, \S 2.2 and \cite{AK}, B.1]
\label{pro AFR-GQA}
(a) For $\xi \in I$, 
we can define a $\Uq$-representation structure on $L_q(\La_{\xi})$ 
by the following formula: For $i \in I$ and $J \in \JJ_{\xi}$, 
\begin{eqnarray}
&& E_iJ=
\begin{cases}
(J\bs\{i+1\}) \sqcup \{i\}, & 
\rm{if $i+1 \in J$ and $i \not\in J$},\\
0, & otherwise, 
\end{cases}
\label{pro AFR-GCA-E} \\
&& F_iJ=
\begin{cases}
(J\bs\{i\}) \sqcup \{i+1\}, & 
\rm{if $i \in J$ and $i+1 \not\in J$},\\
0, & otherwise, 
\label{pro AFR-GCA-F} 
\end{cases}\\
&& K_iJ
=q^{\de(i \in J)-\de(i+1 \in J)}J,
\label{pro AFR-GCA-K} 
\end{eqnarray}
where, for a statement $\th$,
\begin{eqnarray}
 \de(\th):=
\begin{cases}
1, & \rm{if $\th$ is true}, \\
0, & \rm{if $\th$ is false}.
\end{cases}
\end{eqnarray}

(b) $L_q(\La_{\xi})$ is isomorphic to $\Vq(\La_\xi)$ 
as a representation of $\Uq$.
\end{pro}

For $\xi \in I$, we define
\begin{eqnarray}
 J_H^{\xi}:=\{1, 2, \cdots, \xi\}.
\label{def JH}
\end{eqnarray}
Obviously, $J_H^{\xi}$ is a highest-weight vector 
in $L_q(\La_{\xi})$.
\subsection{The evaluation representations of $\tUq$}
In this subsection, 
we introduce the evaluation representations of $\tUq$ 
to consider the fundamental representations of $\tUq$ 
in the next subsection.
\begin{df}
\label{def EQA}
The \it{extended quantum algebra} $U_q^{'}:=U_q^{'}(\sl_{n+1})$ 
is an associative $\bbC(q)$-algebra generated by 
$\{E_{i}, F_{i}, K_{\mu}^{'} \, | \, i \in I, \mu \in P \}$  
with the defining relations    
\begin{eqnarray*}
&& K_{\mu}^{'} K_{\nu}^{'}=K_{\mu+\nu}^{'}, \q K_0^{'}=1,
\q K_{\mu}^{'} E_j K_{-\mu}^{'}
=q^{\langle \mu, \al_i^{\vee} \rangle}E_j,  
\q K_{\mu}^{'} F_j K_{-\mu}^{'}
=q^{-\langle \mu, \al_i^{\vee} \rangle}F_j, \\ 
&& E_{i} F_{j}- F_{j} E_{i} = \de_{i,j} 
\frac{K_{\al_i}^{'}-K_{-\al_i}^{'}}{q-q^{-1}}, \\
&& \sum_{r=0}^{1-\fra_{i,j}} (-1)^r E_{i}^{(r)} E_{j}
E_{i}^{(1-\fra_{i,j}-r)}=
\sum_{r=0}^{1-\fra_{i,j}} (-1)^r F_{i}^{(r)} F_{j} 
F_{i}^{(1-\fra_{i,j}-r)}=0  \q i \neq j, 
\end{eqnarray*}
for $i, j \in I$, $\mu, \nu \in P$.
\end{df}
\begin{rem}
\label{rem EQA}
Let $V$ be a representation of $\Uq$. 
For $i \in I$ and $\mu=\sum_{k \in I}\mu_k\La_k \in P$, we define
\begin{eqnarray*}
 c_{i,\mu}
 :=\frac{1}{n+1}\{(n-i+1)\sum_{k=1}^i k\mu_k
   +i\sum_{k=i+1}^n (n-k+1)\mu_k\},
\label{def cimu}
\end{eqnarray*}
(see (\ref{def fundamental weights})).
We can regard $V$ as a representation of $\Uq^{'}$
by using the following formula: for $v \in V_{\mu}$,  
\begin{eqnarray*}
 K_{\La_i}^{'}v:=q^{c_{i,\mu}}v.
\label{def Kmu}
\end{eqnarray*}
\end{rem}

Now, for $X \in \{E, F\}$, we define 
\begin{eqnarray*}
 X_{\th}^{+}:=[X_n, [X_{n-1}, \cdots, [
 X_2, X_1]_{q^{-1}} \cdots ]_{q^{-1}}, 
 \q X_{\th}^{-}:=[X_1, [X_2, \cdots, [
 X_{n-1}, X_n]_{q^{-1}} \cdots ]_{q^{-1}}, 
\label{def Xth}
\end{eqnarray*}
where $[u,v]_{q^{\pm 1}}=uv-q^{\pm 1}vu$ 
for any $u, v \in U^{'}_q$.
\begin{pro}[\cite{Ji}, \S2 and \cite{CP94a}, Proposition 3.4]
\label{pro EH-GQLA}
For any $\mbf{a} \in \bbC(q)^{\times}$, 
there exist $\bbC(q)$-algebra homomorphisms 
$\mrm{ev}_{\mbf{a}}^{+}$, 
$\mrm{ev}_{\mbf{a}}^{-}:\tUq \arr \Uq^{'}$ such that 
\begin{eqnarray*}
&&\mrm{ev}_{\mbf{a}}^{\pm}(E_i)=E_i, 
\q \mrm{ev}_{\mbf{a}}^{\pm}(F_i)=F_i,
\q \mrm{ev}_{\mbf{a}}^{\pm}(K_i)=K_{\al_i}^{'},  
\q \rm{for} \q i \in I, \\
&&\mrm{ev}_{\mbf{a}}^{\pm}(E_0)
=\mbf{a}K_{\pm (\La_1-\La_n)}^{'}F_{\th}^{\pm}, 
\q \mrm{ev}_{\mbf{a}}^{\pm}(F_0)=(-q)^{n-1}{\mbf{a}}^{-1}
K_{\mp (\La_1-\La_n)}^{'}E_{\th}^{\pm}. 
\end{eqnarray*}
\end{pro}
From Remark \ref{rem EQA} and Proposition \ref{pro EH-GQLA}, 
we can regard any representation of $\Uq$ as a representation of $\tUq$. 
\begin{df} 
\label{def ER-GQLA}
For $\la \in P_{+}$ and $\mbf{a} \in \bbC(q)^{\times}$, we set 
\begin{eqnarray*}
 \mbf{a}_{\la}^{+}
 :=\mbf{a}q^{-c_{1,\la}+c_{n,\la}+n}, \q
 \mbf{a}_{\la}^{-}
 :=(-1)^{n+1}\mbf{a}q^{c_{1,\la}-c_{n,\la}+2n+1}.
\label{def al+-}
\end{eqnarray*} 
We regard $\Vq(\la)$ as a representation of $\tUq$ by using 
$\mrm{ev}_{\mbf{a}_{\la}^{\pm}}^{\pm}$ 
and denote them by $\Vq(\la)_{\mbf{a}}^{\pm}$ 
which are called \it{evaluation representations} of $\Vq(\la)$.
\end{df}
For $i \in I$, $\la \in P_{+}$, and $\bf{a} \in \bbC(q)^{\times}$, 
we define $\mbf{\pi}_i^{\mbf{a},\la}
=(\pi_{i,j}^{\mbf{a},\la}(t))_{j \in I} \in (\bbC(q)_0[t])^n$ whereby 
\begin{eqnarray}
 \pi_{i,j}^{\mbf{a},\la}(t):=
\begin{cases}
\prod_{k=1}^{\la}(1-\mbf{a}q^{\la-2k+1}t), & 
\q \rm{if} \q j=i, \\
1, & \q \rm{if} \q j \neq i,
\end{cases}
\label{def KRP}
\end{eqnarray}
where $\prod_{k=1}^{\la}(1-\mbf{a}q^{\la-2k+1}t):=1$ if $\la=0$.
For $\mbf{\pi}=(\pi_i(t))_{i \in I}$, 
$\mbf{\pi}^{'}=(\pi_i^{'}(t))_{i \in I} \in (\bbC(q)_0[t])^n$, 
we define 
$\mbf{\pi}\mbf{\pi}^{'}:=(\pi_i(t)\pi_i^{'}(t))_{i \in I}$.
For $i \in I$ and $\la=\sum_{i \in I} \la_i\La_i \in P_{+}$, 
we define
\begin{eqnarray}
 \la[i]:=-\sum_{k=1}^{i-1}\la_k+\sum_{k=i+1}^n\la_k-i.
\label{def la[i]}
\end{eqnarray}
The following theorem was proved 
by Chari and Pressley in \cite{CP94a}
(see also \cite{AN}, IV).
\begin{thm}[\cite{CP94a}] 
\label{thm DP-ER-GQLA}
For $\la=\sum_{i \in I}\la_i\La_i \in P_{+}$ 
and $\mbf{a} \in \bbC(q)^{\times}$, 
as representations of $\tUq$, 
\begin{eqnarray*}
 \Vq(\la)_{\mbf{a}}^{\pm} \cong 
 \tVq(\prod_{i \in I}\mbf{\pi}_i^{\mbf{a}q^{\pm\la[i]}, \la_i}). 
\end{eqnarray*}
In particular, we have 
\begin{eqnarray}
 H_{i,1}v_{\la}=\mbf{a}q^{\pm\la[i]-1}[\la_i]_qv_{\la},
\label{fac DP-ER-GQLA}
\end{eqnarray}
where $v_{\la}$ is a highest-weight vector in $\Vq(\la)$.
\end{thm}
\subsection{The fundamental representations of $\tUq$}
For $\xi \in I$ and $\bf{a} \in \bbC(q)^{\times}$, 
we set $\mbf{\pi}_{\xi}^{\mbf{a}}
:=\mbf{\pi}_{\xi}^{\mbf{a}, 1}$. 
Hence $\mbf{\pi}_{\xi}^{\mbf{a}}
=(\pi_{\xi,j}^{\mbf{a}}(t))_{j \in I}$ 
is an element in $(\bbC(q)_0[t])^n$ such that 
\begin{eqnarray}
 \pi_{\xi,j}^{\mbf{a}}(t)=
\begin{cases}
1-\mbf{a}t, & \q \rm{if} \q j=\xi, \\
1, & \q \rm{if} \q j \neq \xi.
\end{cases}
\label{def FP}
\end{eqnarray}
We call $\tVq(\mbf{\pi}_{\xi}^{\mbf{a}})$ 
a \it{fundamental representation} of $\tUq$.

For $\la \in P_{+}$ and $\mbf{a} \in \bbC(q)^{\times}$, 
let $\Vq(\la)_{\mbf{a}}^{\pm}$ be as in Definition \ref{def ER-GQLA}.
For $\xi \in I$, we set 
\begin{eqnarray}
\Vq(\La_{\xi})_{\mbf{a}}
:=\Vq(\La_{\xi})_{\mbf{a}q^{\xi}}^{+}. 
\end{eqnarray}
From Theorem \ref{thm DP-ER-GQLA}, we have
\begin{eqnarray}
\Vq(\La_{\xi})_{\mbf{a}}
=\Vq(\La_{\xi})_{\mbf{a}q^{\xi}}^{+} 
\cong \Vq(\La_{\xi})_{\mbf{a}q^{-\xi}}^{-} 
\cong \tVq(\mbf{\pi}_{\xi}^{\mbf{a}}) 
\q \rm{as representations of $\tUq$}.  
\end{eqnarray}
So we can regard $\Vq(\La_{\xi})_{\mbf{a}}$ 
as the fundamental representation of $\tUq$. 
From Proposition \ref{pro DR+IR}, 
we obtain the following proposition.
\begin{pro} 
\label{pro DFR+IFR-GQLA}
Let $\xi \in I$ and $\mbf{a} \in \bbC(q)^{\times}$. 

(a) There exists an integer $c \in \bbZ$ depending only on $\sl_{n+1}$ 
such that 
$\Vq(\La_{\xi})_{\mbf{a}}^{*}$ is isomorphic to 
$\Vq(\La_{n-\xi+1})_{q^{c} \mbf{a}}$ 
as a representation of $\tUq$. 

(b) There exists a nonzero complex number $\kappa \in \bbC^{\times}$ 
depending only on $\sl_{n+1}$ such that
$\Vq(\La_{\xi})_{\mbf{a}}^{\Om}$ is isomorphic to 
$\Vq(\La_{n-\xi+1})_{q^2\kappa \mbf{a}^{-1}}$ 
as a representation of $\tUq$. 
In particular, $z_{\om_{1,n}\La_{\xi}}$ 
is a pseudo-highest weight vector in $\Vq(\La_{\xi})_{\mbf{a}}^{\Om}$.
\end{pro}

From (\ref{fac DP-ER-GQLA}), for $i \in I$, we have
\begin{eqnarray}
 H_{i,1}z_{\La_{\xi}}
=\La_{\mbf{\pi}_{\xi}^{\mbf{a}}}(H_{i,1})z_{\La_{\xi}}
=\mbf{a}q^{-1} \de_{i, \xi}z_{\La_{\xi}} 
\q \rm{in} \q V_q(\La_{\xi})_{\mbf{a}}.
\label{pro Hs1}
\end{eqnarray}
\begin{pro} 
\label{pro extremal}
Let $\xi \in I$ and $\mbf{a} \in \bbC(q)^{\times}$. 
For $i, j \in I$ such that $i \leq j$, 
let $\om_{i,j}$ be as in (\ref{def om}) 
and let $z_{\La_{\xi}}(\om_{i,j})$ be the 
extremal vector in $V_q(\La_{\xi})_{\mbf{a}}$. 

(a) We have $z_{\La_{\xi}}(\om_{i,j})=z_{\om_{i,j}\La_{\xi}}$.

(b) We have 
\begin{eqnarray*}
&& H_{i-1,1}z_{\La_{\xi}}(\om_{i,j})
=\mbf{a}q^{2j-i-\xi}\de(j-i+2 \leq \xi \leq j)
z_{\La_{\xi}}(\om_{i,j}), \q \rm{if $i \neq 1$},\\
&& H_{j+1,1}z_{\La_{\xi}}(\om_{1,j})
=\mbf{a}q^{j-\xi}\de(\xi \leq j+1)
z_{\La_{\xi}}(\om_{1,j}).
\end{eqnarray*}
\end{pro}

Proof. 
(b) follows from Proposition \ref{pro EV} and (\ref{pro Hs1}). 
So we shall prove (a). 
We obtain 
\begin{eqnarray*}
z_{\La_{\xi}}(\om_{1,1})
=F_1^{\de_{\xi,1}}z_{\La_{\xi}}
=\de_{\xi,1}z_{\La_{\xi}-\al_1}+(1-\de_{\xi,1})z_{\La_{\xi}}
=\de_{\xi,1}z_{s_1\La_{\xi}}+(1-\de_{\xi,1})z_{s_1\La_{\xi}}
=z_{\om_{1,1}\La_{\xi}}.
\end{eqnarray*}
Now we assume $z_{\La_{\xi}}(\om_{i,j})=z_{\om_{i,j}\La_{\xi}}$.
From (\ref{fac vomn1}), 
if $i \neq 1$, we obtain 
\begin{eqnarray*}
z_{\La_{\xi}}(\om_{i-1,j})
=F_{i-1}^{\de(j-i+2 \leq \xi \leq j)}z_{\La_{\xi}}(\om_{i,j})
=F_{i-1}^{\de(j-i+2 \leq \xi \leq j)}z_{\om_{i,j}\La_{\xi}}. 
\end{eqnarray*}
Thus, from (\ref{fac swn1}) and (\ref{pro FR-GCA-F}), we have
\begin{eqnarray*}
 z_{\La_{\xi}}(\om_{i-1,j})
&=&\de(j-i+2 \leq \xi \leq j)z_{\om_{i,j}\La_{\xi}-\al_{i-1}}
+(1-\de(j-i+2 \leq \xi \leq j))z_{\om_{i,j}\La_{\xi}} \\
&=&z_{s_{i-1}\om_{i,j}\La_{\xi}}=z_{\om_{i-1,j}\La_{\xi}}.
\end{eqnarray*}
Similarly, if $i=1$, by using (\ref{fac sw1}) and (\ref{fac vom1}), 
we obtain $z_{\La_{\xi}}(\om_{1,j})=z_{\om_{1,j}\La_{\xi}}$. 
\qed \\

Let $L_q(\La_{\xi})$ and $J_H^{\xi}$ be as in \S 4.2. 
From Proposition \ref{pro AFR-GQA}, Proposition \ref{pro EH-GQLA}, 
and Definition \ref{def ER-GQLA}, we have
\begin{eqnarray}
E_0J_H^{\xi}
=\mbf{a}q^{n+1}(-1)^{\xi-1}\{2, 3, \cdots, \xi, n+1\} 
\q \rm{in} \q L_q(\La_{\xi})_{\mbf{a}} 
\, (\cong \Vq(\La_{\xi})_{\mbf{a}}).
\label{pro E0}
\end{eqnarray}
Moreover, for $J \in \JJ_{\xi}$, we have
\begin{eqnarray*}
&& E_0J=\mbf{a}q^{n-1}(-1)^{\xi}((J\bs\{1\}) \sqcup\{n+1\}), \\
&& F_0J=\mbf{a}^{-1}q^{-n+1}(-1)^{\xi}((J\bs\{n+1\}) \sqcup\{1\}), 
\end{eqnarray*} 
(see \cite{DO} and \cite{AK}). 
\begin{pro} 
\label{pro extremal L}
For $i, j \in I$ such that $i \leq j$,
 let $J_H^{\xi}(\om_{i,j})$ be the 
extremal vector in $L_q(\La_{\xi})_{\mbf{a}}$.
We have  
\begin{eqnarray*}
J_H^{\xi}(\om_{i,j})=
\begin{cases}
J_H^{\xi},& \rm{if $j < \xi$},  \\
\{j+2-\xi, j+3-\xi, \cdots, j+1\}, & 
\rm{if $\xi \leq j$ and $i \leq j+1-\xi$}, \\
\{j+1-\xi, \cdots, i-1, i+1, \cdots, j+1\}, & 
\rm{if $\xi \leq j$ and $j+1-\xi < i$}. 
\end{cases}
\end{eqnarray*}
\end{pro}

Proof. 
This proposition follows from Proposition \ref{pro AFR-GQA} 
and the definition of the extremal vectors in \S 3.4. 
\qed \\

From Lemma 4.2 in \cite{C}, 
we obtain the following lemma.
\begin{lem}
\label{lem HWR+IR=highest}
Let $\xi \in I$, $\mbf{a} \in \bbC(q)^{\times}$, 
and $\mbf{\pi} \in (\bbC(q)_0[t])^n$. 
Let $V$ be a pseudo-highest weight representation of $\tUq$ 
with the pseudo-highest weight $\mbf{\pi}$ 
generated by a pseudo-highest weight vector $v_{\pi}$. 
If $z_{\om_{1,n}\La_{\xi}} \otimes v_{\mbf{\pi}} 
\in \tUq(z_{\La_{\xi}} \otimes v_{\mbf{\pi}})$, 
then $\Vq(\La_{\xi})_{\mbf{a}} \otimes V$ 
is a pseudo-highest weight representation of $\tUq$ 
with the pseudo-highest weight 
$\mbf{\pi}_{\xi}^{\mbf{a}}\mbf{\pi}$ 
generated by a pseudo-highest weight vector 
$z_{\La_{\xi}} \otimes v_{\mbf{\pi}}$. 
\end{lem}

Proof. 
From Lemma \ref{pro tensor highest}, 
it is enough to prove that 
\begin{eqnarray*}
  \Vq(\La_{\xi})_{\mbf{a}} \otimes V
  =\tUq(z_{\La_{\xi}} \otimes v_{\mbf{\pi}}). 
\end{eqnarray*}
Moreover, from the assumption of this Lemma, 
it is enough to prove that
\begin{eqnarray*}
 \Vq(\La_{\xi})_{\mbf{a}} \otimes V
 \subset \tUq(z_{\om_{1,n}\La_{\xi}} \otimes v_{\mbf{\pi}}). 
\end{eqnarray*}
Since $v_{\mbf{\pi}}$ is a pseudo-highest weight vector, 
for any $i_1, i_2, \cdots, i_r \in I$, we have
\begin{eqnarray*}
 E_{i_1}E_{i_2} \cdots E_{i_r}(z_{\om_{1,n}\La_{\xi}} \otimes v_{\mbf{\pi}})
=(E_{i_1}E_{i_2} \cdots E_{i_r}z_{\om_{1,n}\La_{\xi}}) \otimes v_{\mbf{\pi}} 
\q \rm{in} \q \Vq(\La_{\xi})_{\mbf{a}} \otimes V \\
=(F_{i_1}F_{i_2} \cdots F_{i_r}z_{\om_{1,n}\La_{\xi}}) \otimes v_{\mbf{\pi}} 
\q \rm{in} \q \Vq(\La_{\xi})_{\mbf{a}}^{\Om} \otimes V. 
\end{eqnarray*} 
Hence, from Proposition \ref{pro DFR+IFR-GQLA} (b), 
we obtain 
\begin{eqnarray*}
 \tUq(z_{\om_{1,n}\La_{\xi}} \otimes v_{\mbf{\pi}})
\supset (\Uq^{-}z_{\om_{1,n}\La_{\xi}}) \otimes v_{\mbf{\pi}}
=\Vq(\La_{\xi})_{\mbf{a}} \otimes v_{\mbf{\pi}}.
\end{eqnarray*}
Let $Y_{j_1}, Y_{j_2}, \cdots, Y_{j_s} 
\in \{E_i, F_i \, | \, i \in \wt{I}\}$. 
We assume 
$\Vq(\La_{\xi})_{\mbf{a}} \otimes (Y_{j_1} \cdots Y_{j_s}v_{\mbf{\pi}}) 
\subset \tUq(z_{\om_{1,n}\La_{\xi}} \otimes v_{\mbf{\pi}})$. 
Then, for any $\mu \in \WW \La_{\xi}$, we have 
\begin{eqnarray*}
&&q^{(\mu, \al_i)}(z_{\mu} \otimes E_i(Y_{j_1} \cdots Y_{j_s}v_{\mbf{\pi}}))
=E_i(z_{\mu} \otimes (Y_{j_1} \cdots Y_{j_s}v_{\mbf{\pi}}))
-(E_iz_{\mu}) \otimes (Y_{j_1} \cdots Y_{j_s}v_{\mbf{\pi}}) 
\in \tUq(z_{\om_{1,n}\La_{\xi}} \otimes v_{\mbf{\pi}}),\\
&&z_{\mu} \otimes F_i(Y_{j_1} \cdots Y_{j_s}v_{\mbf{\pi}})
=F_i(z_{\mu} \otimes (Y_{j_1} \cdots Y_{j_s}v_{\mbf{\pi}}))
-q^{-(\mu, \al_i)}(F_iz_{\mu}) \otimes (Y_{j_1} \cdots Y_{j_s}v_{\mbf{\pi}}) 
\in \tUq(z_{\om_{1,n}\La_{\xi}} \otimes v_{\mbf{\pi}}). 
\end{eqnarray*}
Therefore we obtain 
\begin{eqnarray*}
 \tUq(z_{\om_{1,n}\La_{\xi}} \otimes v_{\mbf{\pi}}) 
\supset \Vq(\La_{\xi})_{\mbf{a}} \otimes \tUq v_{\mbf{\pi}}
=\Vq(\La_{\xi})_{\mbf{a}} \otimes V.
\end{eqnarray*}
\qed 
\subsection{The fundamental representations of $\Uq(\wt{\sl}_2)$}

We denote the generators 
$X_{1,r}^{\pm}, H_{1,s}, K_1^{\pm 1}$ in $\Uq(\wt{\sl}_2)$ 
(resp. $E_1, F_1, K_1^{\pm 1}$ in $\Uq(\sl_2)$) 
by $X_r^{\pm}$, $H_s$, $K^{\pm 1}$ (resp. $E$, $F$, $K^{\pm 1}$) 
and the fundamental weight $\La_1$ by $\La$.
For $\mbf{a} \in \bbC(q)^{\times}$, 
we denote the $\Uq(\wt{\sl}_2)$-representation 
$\tVq(\mbf{\pi}_1^{\mbf{a}})$ 
(resp. the $\Uq(\sl_2)$-representation $\Vq(\La)$) 
by $\tVq(\mbf{\pi}^{\mbf{a}})$ (resp. $\Vq(1)$). 
Then $\Vq(1)_{\mbf{a}}$ ($\cong \tVq(\mbf{\pi}^{\mbf{a}})$) 
has the following realization: 
\begin{eqnarray}
\Vq(1)_{\mbf{a}}=\bbC(q)z_{\La} \oplus \bbC(q)z_{-\La},
\q X_r^{-}z_{\La}=\mbf{a}^rz_{-\La}, 
\q X_r^{+}z_{-\La}=\mbf{a}^rz_{\La}, 
\q X_r^{+}z_{\La}=X_r^{-}z_{-\La}=0, 
\label{def sl2X}
\end{eqnarray}
for any $r \in \bbZ$.
Moreover, for $m \in \bbN$, we have 
\begin{eqnarray}
\Psi_m^{+}z_{\La}=\mbf{a}^m(q-q^{-1})z_{\La}, 
\q \Psi_0^{+}z_{\La}=Kz_{\La}=qz_{\La}.
\label{def sl2Psi}
\end{eqnarray}
\section{Tensor product of the fundamental representations for 
the quantum loop algebras: the generic case}
\setcounter{equation}{0}
\renewcommand{\theequation}{\thesection.\arabic{equation}}
\subsection{Irreducibility: the $\tUq$ case}
In this subsection, we review the results and proofs in \cite{C} 
that will be needed later. 

For $i \in I$, 
let $\tUq^{(i)}$ (resp. $\Uq^{(i)}$) 
be the $\bbC(q)$-subalgebra of $\tUq$ (resp. $\Uq$) 
generated by $\{X_{i,r}^{\pm}, H_{i,s}, K_i^{\pm 1} \, | 
\, r \in \bbZ, s \in \bbZ^{\times}\}$ 
(resp. $\{E_i, F_i, K_i^{\pm 1}\}$). 
There exist $\bbC(q)$-algebra homomorphisms 
$\wt{\iota}: \Uq(\wt{\sl}_2) \arr \tUq^{(i)}$ 
and $\iota: \Uq(\sl_2) \arr \Uq^{(i)}$ such that 
\begin{eqnarray*}
&& \wt{\iota}(X_r^{\pm})=X_{i,r}^{\pm}, \q 
\wt{\iota}(H_s)=H_{i,s}, \q 
\wt{\iota}(K^{\pm 1})=K_i^{\pm 1}, \\
&& \iota(E)=E_i, \q \iota(F)=F_i, \q \iota(K^{\pm 1})=K_i^{\pm 1},
\end{eqnarray*}
for any $i \in I$, $r \in \bbZ$, and $s \in \bbZ^{\times}$.
Hence, for any $\tUq^{(i)}$-representation 
(resp. $\Uq^{(i)}$-representation) $V$, 
we can regard $V$ as a $\Uq(\wt{\sl}_2)$-representation 
(resp. $\Uq(\sl_2)$-representation).
\begin{lem}
\label{lem TP-FR-GQLA}
Let $\xi \in I$ and $\mbf{a} \in \bbC(q)^{\times}$. 
For any $i, j \in I$ such that $i \leq j$,
as representations of $\Uq(\wt{\sl}_2)$, 
\begin{eqnarray*}
&&\tUq^{(i-1)}z_{\om_{i,j}\La_{\xi}} \cong
  \begin{cases}
    \Vq(1)_{\mbf{a}q^{2j-\xi-i+1}},
    & \rm{if} \q j-i+2 \leq \xi \leq j, \\
    \Vq(0)_{\mbf{a}}, & otherwise, 
  \end{cases}
  \q \rm{if} \q i \neq 1, \\
&&\tUq^{(j+1)}z_{\om_{1,j}\La_{\xi}} \cong
 \begin{cases}
   \Vq(1)_{\mbf{a}q^{j-\xi+1}},
   & \rm{if} \q 1 \leq \xi \leq j+1, \\
   \Vq(0)_{\mbf{a}}, & otherwise. 
 \end{cases}
\end{eqnarray*}
\end{lem}

Proof. 
We assume $i \neq 1$. 
We can prove the case of $i=1$ similarly.
We set $\mu:=\om_{i,j}\La_{\xi}$. 
From (\ref{fac EV}), for any $r \in \bbZ$, 
$X_{i-1, r}^{+}z_{\mu}=0$. 
Hence we have
\begin{eqnarray*}
  \tUq^{(i-1)}z_{\mu} 
  &\subset& (\sum_{r_1, \cdots, r_m \in \bbZ, m \in \bbN} 
  \bbC(q)X_{i-1, r_1}^{-} \cdots X_{i-1, r_m}^{-}z_{\mu})
  \oplus \bbC(q)z_{\mu} \\
  &\subset& \bigoplus_{m \in \bbN}
  (\Vq(\La_{\xi})_{\mu-m\al_{i-1}})
  \oplus \bbC(q)z_{\mu}.
\end{eqnarray*}
Since $\La_{\xi}$ is a minuscule weight (see \S 4.1), 
for any $\nu \in \WW \La_{\xi}$, we have
$| \langle \nu, \al_{i-1}^{\vee} \rangle | \leq 1$. 
Thus we obtain
\begin{eqnarray*}
 | \langle \mu-m\al_{i-1}, \al_{i-1}^{\vee} \rangle |
 =| \langle \mu, \al_{i-1}^{\vee} \rangle -2m|
 \geq 2m-1 \q \rm{for} \q m \in \bbN.
\end{eqnarray*}
So if $m \geq 2$, $\Vq(\La_{\xi})_{\mu-m\al_i}=0$.
Hence we obtain
\begin{eqnarray*}
  \tUq^{(i-1)}z_{\mu} 
  \subset \Vq(\La_{\xi})_{\mu-\al_{i-1}} \oplus \bbC(q)z_{\mu}.
\end{eqnarray*}
Since $\mrm{dim}_{\bbC(q)}(\Vq(\La_{\xi})_{\mu-\al_{i-1}}) \leq 1$ 
from (\ref{pro dim of FR}), we have 
$\mrm{dim}_{\bbC(q)}(\tUq^{(i-1)}z_{\mu}) \leq 2$.
So $\tUq^{(i-1)}z_{\mu} \cong \Vq(0)_{\mbf{a}}$ 
or there exists an element $\mbf{b} \in \bbC(q)^{\times}$ such that 
$\tUq^{(i-1)}z_{\mu} \cong \Vq(1)_{\mbf{b}}$ 
as representations of $\Uq(\wt{\sl}_2)$. 
From (\ref{fac swn1}), we have 
\begin{eqnarray*}
 \langle \mu, \al_{i-1}^{\vee} \rangle
 =\langle \om_{i,j}\La_{\xi}, \al_{i-1}^{\vee} \rangle
 =\de(j-i+2 \leq \xi \leq j). 
\end{eqnarray*}
Therefore, from (\ref{pro FR-GCA-F}), we obtain 
\begin{eqnarray*}
 \tUq^{(i-1)}z_{\mu} \cong
 \begin{cases}
  \Vq(1)_{\mbf{b}},
  & \rm{if} \q j-i+2 \leq \xi \leq j, \\
  \Vq(0)_{\mbf{a}}, & otherwise. 
 \end{cases}
\end{eqnarray*}
We assume $j-i+2 \leq \xi \leq j$. 
From (\ref{pro Hs1}), we have 
\begin{eqnarray*}
 H_{i-1,1}z_{\mu}=\mbf{b}q^{-1}z_{\mu} 
\q \rm{in} \q \Vq(1)_{\mbf{b}}.
\end{eqnarray*}
On the other hand, from Proposition \ref{pro extremal} (b), we have
\begin{eqnarray*}
H_{i-1,1}z_{\mu}
=\mbf{a}q^{2j-i-\xi}z_{\mu} \q \rm{in} \q \Vq(\La_{\xi})_{\mbf{a}}.
\end{eqnarray*}
Therefore we obtain $\mbf{b}=\mbf{a}q^{2j-i-\xi+1}$.
\qed \\

Let $m \in \bbN$ ($m \geq 2$) and 
$\mbf{a}=(\mbf{a}_1, \cdots, \mbf{a}_m) \in (\bbC(q)^{\times})^m$. 
For $1 \leq r, s \leq m$, we define 
\begin{eqnarray}
  A_{r,s}^q(\mbf{a}):=q^{m-s}\mbf{a}_s^r+
  \sum_{k=1}^{r-1}\mbf{a}_s^{r-k} d_{k,s}^q(\mbf{a}), 
  \q A^q(\mbf{a}):=(A_{r,s}^q(\mbf{a}))_{r,s=1}^m, 
\label{def A}
\end{eqnarray}
where, around $u=0$, 
\begin{eqnarray}
  \sum_{k=0}^{\infty} d_{k,s}^q(\mbf{a})u^k
  :=q^m \frac{(1-q^{-2}\mbf{a}_{s+1}u) \cdots (1-q^{-2}\mbf{a}_mu)}
  {(1-\mbf{a}_{s+1}u) \cdots (1-\mbf{a}_mu)} 
  \q \rm{in} \q \bbC(q)[[u]],
\label{def d}
\end{eqnarray}
(see (\ref{fac Psi})) 
and $d_{k,m}^q(\mbf{a}):=0$ for any $k \in \bbZ_{+}$. 
So we have
\begin{eqnarray*}
  \sum_{k=0}^{\infty} d_{k,s}^q(\mbf{a})u^k
  =\prod_{p=s+1}^m\{q+(q-q^{-1})
  (\mbf{a}_pu+\mbf{a}_p^2u^2+\mbf{a}_p^3u^3+ \cdots)\}
\end{eqnarray*}
From the proof of Lemma 4.10 in \cite{CP91}, 
we obtain the following lemma. 
\begin{lem}[\cite{CP91}, \S 4.10]
\label{lem detA}
\begin{eqnarray}
  \mrm{det}(A^q(\mbf{a}))
  =(\prod_{k=1}^m \mbf{a}_k)
 (\prod_{1 \leq s< t \leq m}(q^{-1}\mbf{a}_t-q\mbf{a}_s)).   
\end{eqnarray}
\end{lem}
The following theorem is the special case of Theorem 4.4 in \cite{C}.
\begin{thm}[\cite{C}, Theorem 4.4] 
\label{thm TP-FR-GQLA}
Let $m \in \bbN$, $\xi_1, \cdots, \xi_m \in I$, 
and $\mbf{a}_1, \cdots, \mbf{a}_m \in \bbC^{\times}$. 
We assume that for any $1 \leq k < k^{'} \leq m$ 
and $\mrm{max}(\xi_k, \xi_{k^{'}}) \leq t 
\leq \mrm{min}(\xi_k+\xi_{k^{'}}-1, n)$, 
\begin{eqnarray*}
\frac{\mbf{a}_{k^{'}}}{\mbf{a}_k}
\neq q^{2t-\xi_k-\xi_{k^{'}}+2}. 
\end{eqnarray*}
Then $V_q(\La_{\xi_1})_{\mbf{a}_1} \otimes 
\cdots \otimes V_q(\La_{\xi_m})_{\mbf{a}_m}$ 
is a pseudo-highest weight representation of $\tUq$ 
with the pseudo-highest weight 
$\mbf{\pi}^{\mbf{a}_1}_{\xi_1} \cdots \mbf{\pi}^{\mbf{a}_m}_{\xi_m}$ 
generated by a pseudo-highest weight vector 
$z_{\La_{\xi_1}} \otimes \cdots \otimes z_{\La_{\xi_m}}$.
\end{thm}

Proof. 
We prove this theorem by using the method in \cite{C} and \cite{CP91}. 
From Proposition \ref{pro tensor highest}, 
it is enough to prove that 
\begin{eqnarray*}
 V_q(\La_{\xi_1})_{\mbf{a}_1} \otimes 
\cdots \otimes V_q(\La_{\xi_m})_{\mbf{a}_m}
=\tUq (z_{\La_{\xi_1}} \otimes \cdots \otimes z_{\La_{\xi_m}}).
\end{eqnarray*}
We shall prove this claim by the induction on $m$. 
If $m=1$, we have nothing to prove. 
So we assume $m>1$ and the case of $(m-1)$ holds. 
We set 
\begin{eqnarray*}
 V^{'}:=V_q(\La_{\xi_2})_{\mbf{a}_2} \otimes 
\cdots \otimes V_q(\La_{\xi_m})_{\mbf{a}_m}, 
\q z^{'}:=(z_{\La_{\xi_2}} \otimes \cdots \otimes z_{\La_{\xi_m}}).
\end{eqnarray*}
From Proposition \ref{pro tensor highest} 
and the assumption of the induction on $m$, 
$V^{'}$ is a pseudo-highest weight representation of $\tUq$ 
with the pseudo-highest weight 
$\mbf{\pi}^{\mbf{a}_2}_{\xi_2} \cdots \mbf{\pi}^{\mbf{a}_m}_{\xi_m}$ 
generated by a pseudo-highest weight vector $z^{'}$. 
Hence, from Lemma \ref{lem HWR+IR=highest}, 
it is enough to prove that 
\begin{eqnarray*}
z_{\om_{1,n}\La_{\xi_1}} \otimes z^{'} 
\in \tUq(z_{\La_{\xi_1}} \otimes z^{'}).
\end{eqnarray*}
We shall prove that
\begin{eqnarray}
z_{\om_{i,j}\La_{\xi_1}} \otimes z^{'} 
\in \tUq(z_{\La_{\xi_1}} \otimes z^{'}),
\label{thm TP-FR-GQLA claim}
\end{eqnarray} 
for any $i, j \in I$ such that $i \leq j$. 
We set $z_{\om_{1,0}\La_{\xi_1}}:=z_{\La_{\xi_1}}$. 
We define a total order in 
$I^{\leq}:=\{(i,j) \, | \, 1 \leq i \leq j \leq n\} \sqcup \{(1,0)\}$ 
whereby 
\begin{eqnarray}
  (i-1,j)>(i,j), \q \rm{and} \q (j+1,j+1)>(1,j),
\label{def TO2} 
\end{eqnarray}
for $2 \leq i \leq n$ and $0 \leq j \leq n$. 
We shall prove (\ref{thm TP-FR-GQLA claim}) by the induction on $(i,j)$. 
If $(i,j)=(1,0)$, we have nothing to prove. 
So we assume $(i,j) \neq (1,0)$ and the case of $(i,j)$ holds.
We also assume $i \neq 1$. 
We can prove the case of $i=1$ similarly.
From (\ref{fac swn1}), we have 
$\langle \om_{i,j}, \al_{i-1}^{\vee} \rangle
=\de(j-i+2 \leq \xi \leq j)$. 
So if $\xi_1<j-i+2$ or $\xi_1 > j$, 
then $\om_{i-1,j}=s_{i-1}\om_{i,j}=\om_{i,j}$. 
Hence 
\begin{eqnarray*}
z_{\om_{i-1,j}\La_{\xi_1}} \otimes z^{'}
=z_{\om_{i,j}\La_{\xi_1}} \otimes z^{'}
\in \tUq(z_{\La_{\xi_1}} \otimes z^{'}),  
\end{eqnarray*}
by the induction on $(i,j)$. 
So we assume $j-i+2 \leq \xi_1 \leq j$.

\it{Case 1)} 
In the case of $\xi_2=\cdots =\xi_m=i-1$: 
From Proposition \ref{pro C02-2.2} (b), 
for $r \in \bbZ$, we have
\begin{eqnarray}
&&\De_{H}(X_{i-1,r}^{-})(z_{\om_{i,j}\La_{\xi_1}} \otimes   z^{'})
  -z_{\om_{i,j}\La_{\xi_1}} \otimes (\De_{H}(X_{i-1,r}^{-})z^{'}) \no \\ 
&&=(\De_{H}(X_{i-1,r-k}^{-})z_{\om_{i,j}\La_{\xi_1}})
  \otimes (\De_{H}(K_{i-1})z^{'}) 
  +\sum_{k=1}^{r-1}(X_{i-1,r-k}^{-}z_{\om_{i,j}\La_{\xi_1}})
  \otimes (\De_{H}(\Psi_{i-1,k}^{+})z^{'}). 
\label{thm TP-FR-GQLA 4}
\end{eqnarray}
We have 
\begin{eqnarray}
  \De_{H}(K_{i-1})z^{'}=q^{m-1}z^{'}.
\label{thm TP-FR-GQLA 5}
\end{eqnarray}
We set 
\begin{eqnarray*}
  \wt{\mbf{a}}_1:=\mbf{a}_1 q^{2j-\xi_1-i+1}, 
  \q \wt{\mbf{a}}_2:=\mbf{a}_2, \q \cdots, \q \wt{\mbf{a}}_m:=\mbf{a}_m, 
  \q \wt{\mbf{a}}:=(\wt{\mbf{a}}_1, \cdots, \wt{\mbf{a}}_m).
\end{eqnarray*}
From Lemma \ref{lem TP-FR-GQLA},  
$\tUq^{(i)}z_{\om_{i,j}\La_{\xi_1}} \cong 
\Vq(1)_{\wt{\mbf{a}}_1}$ 
as representations of $\Uq(\wt{\sl}_2)$ 
and $z_{\om_{i,j}\La_{\xi_1}}$ is a pseudo-highest weight vector 
in $\tUq^{(i)}z_{\om_{i,j}\La_{\xi_1}}$. 
Hence, from (\ref{def sl2X}), we have 
\begin{eqnarray}
  \De_H(X_{i-1,r}^{-})z_{\om_{i,j}\La_{\xi_1}}
  =\wt{\mbf{a}}_1^rF_{i-1}z_{\om_{i,j}\La_{\xi_1}}
  =\wt{\mbf{a}}_1^rz_{\om_{i-1,j}\La_{\xi_1}}.
\label{thm TP-FR-GQLA 6}
\end{eqnarray}
From (\ref{fac Psi}),
for $1 \leq k \leq r-1$, we obtain 
\begin{eqnarray}
  \De_H(\Psi_{i-1,k})z^{'}=d_{k,1}^q(\wt{\mbf{a}})z^{'}, 
\label{thm TP-FR-GQLA 7}
\end{eqnarray}
where $d_{k,1}^q(\wt{\mbf{a}})$ be as in (\ref{def d}). 
From (\ref{thm TP-FR-GQLA 4})--(\ref{thm TP-FR-GQLA 7}), 
we obtain
\begin{eqnarray*}
&&\De_{H}(X_{i-1,r}^{-})(z_{\om_{i,j}\La_{\xi_1}} \otimes z^{'})
  -z_{\om_{i,j}\La_{\xi_1}} \otimes 
  (\De_{H}(X_{i-1,r}^{-})z^{'}) \no \\
&&=(q^{m-1}  \wt{\mbf{a}}_1^r 
  +\sum_{k=1}^{r-1} \wt{\mbf{a}}_1^{r-k}
  d_{k,1}^q(\wt{\mbf{a}}))(z_{\om_{i-1,j}\La_{\xi_1}} \otimes z^{'}) 
  =A_{r,1}^q(\wt{\mbf{a}})(z_{\om_{i-1,j}\La_{\xi_1}} \otimes z^{'}).
\end{eqnarray*}
By repeating this procedure $m$-times, we obtain
\begin{eqnarray*}
&&\De_{H}(X_{i-1,r}^{-})(z_{\om_{i,j}\La_{\xi_1}} \otimes z^{'})\no \\
&&=A_{r,1}^q(\wt{\mbf{a}})(z_{\om_{i-1,j}\La_{\xi_1}} \otimes z^{'}) 
  +\sum_{s=2}^m A_{r,s}^q(\wt{\mbf{a}})
  (z_{\om_{i,j}\La_{\xi_1}} \otimes 
  z_{\La_{\xi_2}} \otimes \cdots \otimes
  (F_{i-1}z_{\La_{\xi_s}}) \otimes \cdots \otimes z_{\La_{\xi_m}}).
\end{eqnarray*}
Since $\De_{H}(X_{i-1,r}^{-})(z_{\om_{i,j}\La_{\xi_1}} \otimes z^{'})
\in \tUq(z_{\La_{\xi_1}} \otimes z^{'})$ 
from the assumption of the induction on $(i,j)$, we obtain 
\begin{eqnarray*}
  \mrm{det}(A^q(\wt{\mbf{a}}))
  (z_{\om_{i-1,j}\La_{\xi_1}} \otimes z^{'})
  \in \tUq (z_{\La_{\xi_1}} \otimes z^{'}),
\end{eqnarray*}
where $A^q(\wt{\mbf{a}})=(A^q_{r,s}(\wt{\mbf{a}})))_{r,s=1}^m$ 
be as in (\ref{def A}). 
Hence, from Lemma \ref{lem detA} and the assumption of Case 1, 
we have
\begin{eqnarray}
  (\prod_{k=2}^m (\mbf{a}_k-\mbf{a}_1q^{2j-\xi_1-\xi_k+2}))
  (\prod_{2 \leq k < k^{'} \leq m} 
  (\mbf{a}_{k^{'}}-\mbf{a}_kq^2))
  (z_{\om_{i-1,j}\La_{\xi_1}} \otimes z^{'})
  \in \tUq (z_{\La_{\xi_1}} \otimes z^{'}).
\label{thm TP-FR-GQLA 8}
\end{eqnarray}
From the assumption of this theorem, 
for any $1 \leq k \leq m$, we have
\begin{eqnarray}
  \mbf{a}_k-\mbf{a}_1q^{2j-\xi_1-\xi_k+2} \neq 0.
\label{thm TP-FR-GQLA 9}
\end{eqnarray}
Moreover, 
if $\xi_k=\xi_{k^{'}}$, we have
$\mrm{max}(\xi_k, \xi_{k^{'}}) \leq \xi_k=\xi_k^{'} 
\leq \mrm{min}(\xi_k+\xi_{k^{'}}-1, n)$. 
Thus, from the assumption of this theorem and Case 1, 
we have  
\begin{eqnarray}
  \mbf{a}_{k^{'}}-\mbf{a}_kq^2 \neq 0.
\label{thm TP-FR-GQLA 10}
\end{eqnarray}
for any $2 \leq k < k^{'} \leq m$.
Therefore, from (\ref{thm TP-FR-GQLA 8})--(\ref{thm TP-FR-GQLA 10}), 
we obtain 
\begin{eqnarray*}
  z_{\om_{i-1,j}\La_{\xi_1}} \otimes z^{'}
  \in \tUq (z_{\La_{\xi_1}} \otimes z^{'}).
\end{eqnarray*}

\it{Case 2) There exists an integer $m^{'}$ ($2 \leq m^{'} \leq m$) 
such that $\xi_{m^{'}} \neq i-1$}: 
We set 
\begin{eqnarray}
  M:=\{2 \leq m^{'} \leq m \, | \, \xi_{m^{'}} = i-1\}.
\label{thm TP-FR-GQLA 12} 
\end{eqnarray}
If $M=\emptyset$, then $F_{i-1}z^{'}=0$. Hence we obtain
\begin{eqnarray*}
  \tUq (z_{\La_{\xi_1}} \otimes z^{'}) 
  \ni F_{i-1}(z_{\om_{i,j}\La_{\xi_1}} \otimes z^{'})
  =q^{-1} (z_{\om_{i-1,j}\La_{\xi_1}} \otimes z^{'}).
\end{eqnarray*}
We assume $M \neq  \emptyset$. 
For $m^{'} \in \{1, 2, \cdots, m\}$ such that $\xi_{m^{'}} \neq i-1$, 
we have
\begin{eqnarray*}
  \De_H(X_{i-1,r}^{-})
  (z_{\xi_{m^{'}}} \otimes z_{\xi_{m^{'}+1}} 
  \otimes \cdots \otimes z_{\xi_m}) 
  =z_{\xi_{m^{'}}} \otimes \De_H(X_{i-1,r}^{-})(z_{\xi_{m^{'}+1}} 
  \otimes \cdots \otimes z_{\xi_m}),
\end{eqnarray*}
(see (\ref{thm TP-FR-GQLA 4})).
Hence, in a similar way to the proof of Case 1, we obtain
\begin{eqnarray*}
  (\prod_{k \in M} (\mbf{a}_k-\mbf{a}_1q^{2j-\xi_1-\xi_k+2}))
  (\prod_{k, k^{'} \in M, k<k^{'}} 
  (\mbf{a}_{k^{'}}-\mbf{a}_kq^2))
  (z_{\om_{i-1,j}\La_{\xi_1}} \otimes z^{'})
  \in \tUq (z_{\La_{\xi_1}} \otimes z^{'}).
\label{thm TP-FR-GQLA 11}
\end{eqnarray*}
Therefore, from the assumption of this theorem, we obtain 
\begin{eqnarray*}
  z_{\om_{i-1,j}\La_{\xi_1}}\otimes z^{'}
  \in \tUq (z_{\La_{\xi_1}} \otimes z^{'}).
\end{eqnarray*}
\qed 
\begin{rem}
\label{rem irreducibility tUqsl2}
Since the comultiplication $\De_H$ in this paper is 
slightly different from the one in \cite{CP91}, 
$\mrm{det}(A^q(\mbf{a}))$ is different from 
the one in the proof of Lemma 4.10 in \cite{CP91}.
\end{rem}
\begin{cor} 
\label{cor TP-FR-GQLA}
Let $m \in \bbN$, $\xi_1, \cdots, \xi_m \in I$, 
and $\mbf{a}_1, \cdots, \mbf{a}_m \in \bbC(q)^{\times}$. 
We assume that for any $1 \leq k \neq k^{'} \leq m$ 
and $\mrm{max}(\xi_k, \xi_{k^{'}}) \leq t 
\leq \mrm{min}(\xi_k+\xi_{k^{'}}-1, n)$, 
\begin{eqnarray*}
\frac{\mbf{a}_{k^{'}}}{\mbf{a}_k}
\neq q^{\pm (2t-\xi_k-\xi_{k^{'}}+2)}. 
\end{eqnarray*}
Then $V_q(\La_{\xi_1})_{\mbf{a}_1} \otimes 
\cdots \otimes V_q(\La_{\xi_m})_{\mbf{a}_m}$ 
is an irreducible representation of $\tUq$.
\end{cor}

Proof. 
We set 
$V:=V_q(\La_{\xi_1})_{\mbf{a}_1} \otimes 
\cdots \otimes V_q(\La_{\xi_m})_{\mbf{a}_m}$. 
It is enough to prove that $V^{\Om}$ is irreducible. 
Since $(V^{\Om})^{\Om} \cong V$, 
from Theorem \ref{thm TP-FR-GQLA}, 
$(V^{\Om})^{\Om}$ is a pseudo-highest representation of $\tUq$. 
Hence, from Proposition \ref{pro DR+IR=irreducible}, 
it is enough to prove that $(V^{\Om})^{*}$ is 
a pseudo-highest representation of $\tUq$.
From Proposition \ref{pro DFR+IFR-GQLA} (b) and 
(\ref{pro DR+IR tensor}), 
there exists a nonzero complex number $\kappa \in \bbC^{\times}$ 
such that
\begin{eqnarray*}
 V^{\Om} \cong \Vq(\La_{\xi_m})_{q^2\kappa \mbf{a}_m^{-1}} \otimes \cdots 
 \otimes \Vq(\La_{\xi_1})_{q^2\kappa \mbf{a}_1^{-1}}.
\end{eqnarray*}
Thus, From Proposition \ref{pro DFR+IFR-GQLA} (a), 
there exists an integer $c \in \bbZ$ such that
\begin{eqnarray*}
 (V^{\Om})^{*}  \cong \Vq(\La_{\xi_1})_{q^{c+2}\kappa \mbf{a}_1^{-1}} 
 \otimes \cdots 
 \otimes \Vq(\La_{\xi_m})_{q^{c+2}\kappa\mbf{a}_m^{-1}}.
\end{eqnarray*}
From the assumption of this corollary, 
for $1 \leq k < k^{'} \leq m$, we have
\begin{eqnarray*}
 \frac{q^{c+2} \kappa \mbf{a}_{k^{'}}^{-1} }
 {q^{c+2} \kappa \mbf{a}_k^{-1}}
 =\frac{\mbf{a}_k}{\mbf{a}_{k^{'}}}
 \neq q^{2t-\xi_k-\xi_{k^{'}}+2}.
\end{eqnarray*}
Therefore, from Theorem \ref{thm TP-FR-GQLA}, 
$(V^{\Om})^{*}$ is a pseudo-highest representation of $\tUq$. 
\qed
\subsection{The $R$-matrices of the fundamental 
 representations of $\tUq$}
In this subsection, we regard $\Vq(\La_{\xi})$ 
as $L_q(\La_{\xi})$ for $\xi \in I$ (see \S 4.2). 
We review the $R$-matrices of the fundamental 
representations of $\tUq$ introduced in \cite{DO}. 

For $\xi, \ze \in I$, as representations of $\Uq$, 
\begin{eqnarray}
 \Vq(\La_{\xi}) \otimes \Vq(\La_{\ze}) 
 \cong \bigoplus_{k=\mrm{max}(0, \xi+\ze-n-1)}^{\mrm{min}(\xi,\ze)}
 \Vq(\La_{\xi+\ze-k}+\La_k).
\label{fac tensor decomposition}
\end{eqnarray}
For $k \in I$ such that 
$\mrm{max}(0, \xi+\ze-n-1) \leq k \leq \mrm{min}(\xi,\ze)$, 
we define
\begin{eqnarray}
w_k^{(\xi,\ze)}
:=\sum_{J \subset (J^{(\xi+\ze-k)}\bs J^{(k)}), \#(J)=\xi-k} 
(-q)^{\sum_{j \in J}j-\frac{(\xi+k+1)(\xi-k)}{2}}
(J^{(k)} \sqcup J) \otimes (J^{(\xi+\ze-k)}\bs J), 
\label{fac HWV-FR-GQA}
\end{eqnarray}
where $J^{(i)}:=\{1, 2, \cdots, i\}$ $(i \in \wh{I})$ 
and $J^{(0)}:=\emptyset$. 

Let $L_k^{(\xi,\ze)}$ be the $\Uq$-subrepresentation of 
$\Vq(\La_{\xi}) \otimes \Vq(\La_{\ze})$ 
such that $L_k^{(\xi,\ze)} \cong \Vq(\La_{\xi+\ze-k}+\La_k)$. 
Then $w_k^{(\xi,\ze)}$ is a highest-weight vector in $L_k^{(\xi,\ze)}$. 
For $k \in I$ such that
$\mrm{max}(0, \xi+\ze-n-1) \leq k \leq \mrm{min}(\xi,\ze)$, 
there exists a $\Uq$-homomorphism 
$\ol{P}_k: \Vq(\La_{\xi}) \otimes \Vq(\La_{\ze}) 
\arr \Vq(\La_{\ze}) \otimes \Vq(\La_{\xi})$ such that 
\begin{eqnarray}
  \ol{P}_k(w_{k^{'}}^{(\xi,\ze)})
  =\de_{k,k^{'}}w_k^{(\ze,\xi)}.
\label{def Pbar}
\end{eqnarray}
For $\mbf{a}, \mbf{b} \in \bbC(q)^{\times}$, 
we define a $\bbC(q)$-linear map 
$\ol{R}_{\xi,\ze}(\mbf{a}, \mbf{b})
: \Vq(\La_{\xi})_{\mbf{a}} \otimes \Vq(\La_{\ze})_{\mbf{b}} 
\arr \Vq(\La_{\ze})_{\mbf{b}} \otimes \Vq(\La_{\xi})_{\mbf{a}}$ whereby 
\begin{eqnarray}
 \ol{R}_{\xi,\ze}(\mbf{a}, \mbf{b})
 :=\sum_{k=\mrm{max}(0, \xi+\ze-n-1)}^{\mrm{min}(\xi,\ze)}
 \ol{\rho}_k \ol{P}_k,
\label{def olR}
\end{eqnarray} 
where
\begin{eqnarray*}
 \frac{\ol{\rho}_{k-1}}{\ol{\rho}_k}
 :=\frac{\mbf{b}-q^{\xi+\ze-2k+2}\mbf{a}}
 {\mbf{a}-q^{\xi+\ze-2k+2}\mbf{b}}, 
 \q \ol{\rho}_{\mrm{min}(\xi,\ze)}:=1.
\end{eqnarray*}
We call $\ol{R}_{\xi,\ze}(\mbf{a}, \mbf{b})$ a \it{$R$-matrix} 
of $\Vq(\La_{\xi})_{\mbf{a}} \otimes \Vq(\La_{\ze})_{\mbf{b}}$.
From \cite{DO}, \S 2.3 (and (\ref{pro E0}) in this paper), 
we obtain the following theorem. 
\begin{thm}[\cite{DO}, \S 2.3] 
\label{thm R-matrix-FR-GQLA}
$\ol{R}_{\xi,\ze}(\mbf{a}, \mbf{b})$ is a $\tUq$-homomorphism. 
\end{thm}
We have 
\begin{eqnarray*}
 &&\ol{\rho}_k
 =\prod_{p=k+1}^{\mrm{min}(\xi,\ze)}
 \frac{\mbf{b}-q^{\xi+\ze-2k+2}\mbf{a}}
 {\mbf{a}-q^{\xi+\ze-2k+2}\mbf{b}}
 =\prod_{p=1}^{\mrm{min}(\xi,\ze)-k}
 \frac{\mbf{b}-q^{2p+|\xi-\ze|}\mbf{a}}
 {\mbf{a}-q^{2p+|\xi-\ze|}\mbf{b}}, 
 \q (p \mapsto \mrm{min}(\xi, \ze)+1-p), \\
 &&\ol{R}_{\xi,\ze}(\mbf{a}, \mbf{b})
 =\sum_{k=0}^{\mrm{min}(\xi,\ze, n+1-\xi, n+1-\ze)}
 \ol{\rho}_{\mrm{min}(\xi,\ze)-k} \ol{P}_{\mrm{min}(\xi,\ze)-k}, 
 \q (k \mapsto \mrm{min}(\xi, \ze)-k).
\end{eqnarray*}
We set 
\begin{eqnarray}
 && R_{\xi,\ze}(\mbf{a}, \mbf{b})
 :=(\prod_{p=1}^{\mrm{min}(\xi,\ze)}(\mbf{a}-q^{2p+|\xi-\ze|}\mbf{b}))
 \ol{R}_{\xi,\ze}(\mbf{a}, \mbf{b}), \\
 && \rho_k
 :=(\prod_{p=1}^{\mrm{min}(\xi,\ze)}(\mbf{a}-q^{2p+|\xi-\ze|}\mbf{b}))
 \ol{\rho}_{\mrm{min}(\xi, \ze)-k}, 
 \q P_k:=\ol{P}_{\mrm{min}(\xi, \ze)-k}.
 \label{def rhoP}
\end{eqnarray}
Then we have 
\begin{eqnarray}
 R_{\xi,\ze}(\mbf{a}, \mbf{b})
 =\sum_{k=0}^{\mrm{min}(\xi,\ze, n+1-\xi, n+1-\ze)}
 \rho_k P_k, 
 \qq \rho_k
 =(\prod_{p=1}^k (\mbf{b}-q^{2p+|\xi-\ze|}\mbf{a}))
 (\prod_{p=k+1}^{\mrm{min}(\xi, \ze)}
 (\mbf{a}-q^{2p+|\xi-\ze|}\mbf{b})). 
\label{pro R-matrix}
\end{eqnarray}
\subsection{Reducibility: the $\tUq$ case}
\begin{pro} 
\label{pro Reducibility-FR-GQLA}
Let $\xi, \ze \in I$ and $\mbf{a}, \mbf{b} \in \bbC(q)^{\times}$. 
If there exists a 
$1 \leq p_0 \leq \mrm{min}(\xi, \ze, n+1-\xi, n+1-\ze)$ such that 
$\mbf{b}=q^{2p_0+|\xi-\ze|}\mbf{a}$ or $q^{-(2p_0+|\xi-\ze|)}\mbf{a}$, 
then $\Vq(\La_{\xi})_{\mbf{a}} \otimes \Vq(\La_{\ze})_{\mbf{b}}$ 
is a reducible representation of $\tUq$. 
\end{pro}

Proof. 
We assume $\mbf{b}=q^{2p_0+|\xi-\ze|}\mbf{a}$. 
We can prove the case of 
$\mbf{b}=q^{-(2p_0+|\xi-\ze|)}\mbf{a}$ similarly. 
It is enough to prove that 
$R_{\xi,\ze}(\mbf{a}, \mbf{b})$ is neither 
an isomorphism nor a zero map. 
Since $\mbf{b}=q^{2p_0+|\xi-\ze|}\mbf{a}$, 
for $p_0 \leq k \leq \mrm{min}(\xi,\ze,n+1-\xi, n+1-\ze)$, 
we have $\rho_k=0$. 
Hence $R_{\xi,\ze}(\mbf{a}, \mbf{b})$ is not an isomorphism. 

Now we assume $R_{\xi,\ze}(\mbf{a}, \mbf{b})=0$.
Then $\rho_k=0$ for any $0 \leq k \leq p_0-1$. 
Thus, for any $0 \leq k \leq p_0-1$, 
there exists a $0 \leq p^{'} \leq p_0-1$ such that 
$\mbf{b}=q^{2p^{'}+|\xi-\ze|}\mbf{a}$ 
or a $p_0 \leq p^{''} \leq \mrm{min}(\xi, \ze)$ such that 
$\mbf{a}=q^{2p^{''}+|\xi-\ze|}\mbf{b}$. 
Since $\mbf{b}=q^{2p_0+|\xi-\ze|}\mbf{a}$, we obtain 
\begin{eqnarray*}
  p=p^{'}, \q \rm{or} \q p_0+p^{''}+|\xi-\ze|=0. 
\end{eqnarray*} 
However this condition never occurs. 
Therefore $R_{\xi,\ze}(\mbf{a}, \mbf{b})$ is not a zero map.
\qed
\subsection{Main theorem: the $\tUq$ case}
\begin{thm} 
\label{thm Main theorem-GQLA}
Let $m \in \bbN$, $\xi_1, \cdots, \xi_m \in I$, 
and $\mbf{a}_1, \cdots, \mbf{a}_m \in \bbC(q)^{\times}$. 
The following conditions (a) and (b) are equivalent. 

(a) 
$\tVq(\pi_{\xi_1}^{\mbf{a}_1}) \otimes 
\cdots \otimes \tVq(\pi_{\xi_m}^{\mbf{a}_m})$ 
is an irreducible representation of $\tUq$.

(b) 
For any $1 \leq k \neq k^{'} \leq m$ 
and $1 \leq t \leq 
\mrm{min}(\xi_k, \xi_{k^{'}}, n+1-\xi_k, n+1-\xi_{k^{'}})$, 
\begin{eqnarray*}
  \frac{\mbf{a}_{k^{'}}}{\mbf{a}_k}
  \neq q^{\pm (2t+|\xi_k-\xi_{k^{'}}|)}. 
\end{eqnarray*}
\end{thm}

Proof. 
(b) is equivalent to the following condition (b)$^{'}$: 

(b)$^{'}$ 
  For any $1 \leq k \neq k^{'} \leq m$ 
  and $\mrm{max}(\xi_k, \xi_{k^{'}}) \leq t 
  \leq \mrm{min}(\xi_k+\xi_{k^{'}}-1, n)$, 
  \begin{eqnarray*}
    \frac{\mbf{a}_k^{'}}{\mbf{a}_k}
    \neq q^{\pm (2t-\xi_k-\xi_{k^{'}}+2)}. 
  \end{eqnarray*}

Therefore this theorem follows from 
Corollary \ref{cor TP-FR-GQLA} 
and Proposition \ref{pro Reducibility-FR-GQLA}. 
\qed 
\section{Tensor product of the fundamental representations for 
the restricted quantum loop algebras}
\setcounter{equation}{0}
\renewcommand{\theequation}{\thesection.\arabic{equation}}

In the rest of this paper, we fix the following notations. 
Let $l$ be an odd integer greater than $2$, 
let $\vep$ be a primitive $l$-th root of unity, 
and let ${\cal A}:=\bbC[t,t^{-1}]$ be the Laurent polynomial ring. 
We regard $\bbC$ as an ${\cal A}$-algebra by the following formula:
\begin{eqnarray*}
  g(q).c:=g(\vep)c \q \rm{for} \q g(q) \in {\cal A}, c \in \bbC, 
\end{eqnarray*}
and denote it by $\bbC_{\vep}$.  

\subsection{Definition of the restricted quantum loop algebras}
 
For $i \in \wt{I}$ and $m \in \bbN$, 
let $E_i^{(m)}$ and $F_i^{(m)}$ be as in (\ref{def divided power}).
Let $\tUAres$ (resp. $\UAres$) be the $\cal{A}$-subalgebra 
of $\tUq$ (resp. $\Uq$) generated by 
$\{E_i^{(m)}, F_i^{(m)}, K_i^{\pm 1} \, | 
\, i \in \wt{I} \, (\rm{resp. $i \in I$}), m \in \bbN\}$. 
For $i \in I$, $r \in \bbZ$, and $m \in \bbN$, we define   
\begin{eqnarray*}
 \left[
  \begin{array}{r}
    K_i ; r\\
    m  \, \, \,
  \end{array}
  \right]
  :=\prod_{p=1}^m 
  \frac{K_iq^{r-p+1}-K_i^{-1}q^{-r+p-1}}{q^p-q^{-p}}.
\end{eqnarray*}
It is known that 
$ \left[
\begin{array}{r}
  K_i ; r\\
  m \, \,
\end{array}
\right] 
\in \UAres$
(see \cite{CP94b}, \S 9.3A). 
Moreover, we have 
\begin{eqnarray*}
  (X_{i,r}^{\pm})^{(m)}:=\frac{1}{[m]_q!}(X_{i,r}^{\pm})^{m}, \q 
  \frac{1}{[s]_q}H_{i,s}, \q \PP_{i,r} \in \tUAres,
\end{eqnarray*}
for $i \in I$, $r \in \bbZ$, and $s \in \bbZ^{\times}$, 
where $\PP_{i,r}$ be as in (\ref{def PP}) 
(see \cite{CP97}, \S 3.1).
\begin{df}
\label{def RQA}
We define 
\begin{eqnarray*}
  \tUres:=\tUAres \otimes_{\AA} \bbC_{\vep}, 
  \q (\rm{resp.} \q \Ures:=\UAres \otimes_{\AA} \bbC_{\vep}),
\end{eqnarray*}
which is called 
a \it{restricted quantum loop algebra} 
(or \it{quantum loop algebra of Lusztig type}) 
(resp. \it{restricted quantum algebra} 
(or \it{quantum algebra of Lusztig type})) 
(see \cite{L89} and \cite{CP97}).
\end{df}
For $i \in \wt{I}$, $j \in I$, $r \in \bbZ$, $s \in \bbZ^{\times}$, 
and $ m \in \bbN$,
we set 
\begin{eqnarray*}
&&e_i^{(m)}:=E_i^{(m)} \otimes 1, \q f_i^{(m)}:=F_i^{(m)} \otimes 1, \q
  k_i:=K_i \otimes 1, \q
  \left[
  \begin{array}{r}
  k_j ; r\\
  m \, \,
  \end{array}
  \right]
  :=\left[
  \begin{array}{r}
  K_j ; r\\
  m \, \,
  \end{array}
  \right], \\
&&(x_{i,r}^{\pm})^{(m)}:=(X_{i,r}^{\pm})^{(m)} \otimes 1, \q
  h_{i,s}:=H_{i,s} \otimes 1, \q
  \psi_{i,r}^{\pm}=\Psi_{i,r}^{\pm} \otimes 1, \q 
  \frp_{i,r}:=\PP_{i,r} \otimes 1 \q \in \tUres.
\end{eqnarray*}
\subsection{The triangular decompositions of $\Ures$ and $\tUres$}
In a similar way to the proof of Lemma 3.4 in \cite{AN}, 
we can prove the following lemma. 
\begin{lem} 
\label{lem A-basis}
Let $V$ be a free $\AA$-module 
and let $\{v_j\}_{j \in J}$ be a $\AA$-basis of $V$, 
where $J$ is an index set. 
Then $\{v_j \otimes 1\}_{j \in J}$ 
is a $\bbC$-basis of $V \otimes_{\AA} \bbC_{\vep}$.
\end{lem}

Proof. 
It is enough to prove that 
$\{v_j \otimes 1\}_{j \in J}$ is $\bbC$-linearly independent in 
$V \otimes_{\AA} \bbC_{\vep}$. 
For $c_j \in \bbC$ ($j \in J$), we assume 
$\sum_{j \in J}c_j(v_j \otimes 1)=0$. 
Then we have $\sum_{j \in J} c_jv_j \in (q-\vep)V$. 
Since $V$ is generated by $\{v_j\}_{j \in J}$ as a $\AA$-module, 
there exist $d_j \in \AA$ ($j \in J$) such that 
\begin{eqnarray*}
  \sum_{j \in J} c_jv_j=(q-\vep)\sum_{j \in J}d_jv_j.  
\end{eqnarray*}
Since $\{v_j\}_{j \in J}$ is $\AA$-linearly independent in $V$, 
for any $j \in J$, we have $c_j=(q-\vep)d_j$. 
Hence there exist $m \in \bbN$ and $c_{j,k} \in \bbC$ 
($-m \leq k \leq m$) such that
\begin{eqnarray*}
  c_j=(q-\vep)\sum_{k=-m}^m c_{j,k}q^k.
\end{eqnarray*}
Therefore we obtain $c_j=0$ for all $j \in J$.
\qed \\

Let $(\tUres)^{+}$ (resp. $(\tUres)^{-}$, $(\tUres)^{0}$) 
be the $\bbC$-subalgebra of $\tUres$ 
generated by $\{e_i^{(m)} \, | \, i \in \wt{I}, m \in \bbN\}$ 
(resp. $\{f_i^{(m)} \, | \, i \in \wt{I}, m \in \bbN\}$, 
$\{\frp_{i,r}, k_i, 
 \left[
\begin{array}{r}
k_i ; 0\\
l \, \,
\end{array}
\right] \, 
|\, i \in I, r \in \bbZ, m \in \bbN\}$). 
Similarly, 
let $(\Ures)^{+}$ (resp. $(\Ures)^{-}$, $(\Ures)^{0}$) 
be the $\bbC$-subalgebra of $\Ures$ 
generated by $\{e_i^{(m)} \, | \, i \in I, m \in \bbN\}$ 
(resp. $\{f_i^{(m)} \, | \, i \in I, m \in \bbN\}$, 
$\{k_i, 
 \left[
\begin{array}{r}
k_i ; 0\\
l \, \,
\end{array}
\right] \, 
|\, i \in I\}$). 
From \cite{L90a}, 
we obtain the triangular decomposition of $\Ures$, 
that is, 
the multiplication defines an isomorphism of $\bbC$-vector spaces: 
\begin{eqnarray}
  (\Ures)^{-} \otimes (\Ures)^0 \otimes (\Ures)^{+} \wt{\arr} \Ures,
\label{TD-RQA}
\end{eqnarray}
(see also \cite{CP94b}, Proposition 9.3.3 and \cite{CP97}, \S 1).
From \cite{CP97}, \S 6, we obtain 
\begin{eqnarray*}
  \tUAres=(\tUAres)^{-} (\tUAres)^0 (\tUAres)^{+}.
\end{eqnarray*}
Hence from (\ref{fac TD-GQLA}), 
we obtain the triangular decomposition of $\tUAres$: 
\begin{eqnarray}
  (\tUAres)^{-} \otimes (\tUAres)^0 \otimes (\tUAres)^{+} 
  \wt{\arr} \tUAres.
\label{TD-AQLA}
\end{eqnarray}
Therefore, from Lemma \ref{lem A-basis}, 
we obtain the triangular decomposition of $\tUres$: 
\begin{eqnarray}
  (\tUres)^{-} \otimes (\tUres)^0 \otimes (\tUres)^{+} 
  \wt{\arr} \tUres. 
\label{TD-RQLA}
\end{eqnarray}
\subsection{The comultiplication of $\Ures$ and $\tUres$}
For $i \in I$ and $m \in \bbN$, we have  
\begin{eqnarray*}
&& \De_H(E_i^{(m)})=\sum_{p=0}^m 
   q^{p(m-p)}E_i^{(m-p)}K_i^p \otimes E_i^{(p)},
   \q \De_H(F_i^{(m)})=\sum_{p=0}^m 
   q^{p(m-p)}F_i^{(p)} \otimes F_i^{(m-p)}K_i^{-p}, \\
&& \ep_{H}(E_i^{(m)})=\ep_H(F_i^{(m)})=0,  \\
&& S_H(E_i^{(m)})=(-1)^m q^{m(m-1)}K_i^{-m}E_i^{(m)}, 
   \q S_H(F_i^{(m)})=(-1)^m q^{-m(m-1)}F_i^{(m)}K_i^m,  
\end{eqnarray*} 
(see \cite{Ja}, \S 3--4). 
Hence $\tUAres$ (resp. $\UAres$) is a Hopf subalgebra of 
$\tUq$ (resp. $\Uq$). 
In particular, 
we can define Hopf algebra structures on $\tUres$ and $\Ures$. 
The comultiplication $\De_H^{\mrm{res}}$, counit $\ep_H^{\mrm{res}}$, 
and antipode $S_H^{\mrm{res}}$ of $\tUres$ (resp. $\Ures$) are given by 
\begin{eqnarray*}
  \De_H^{\mrm{res}}=\De_H \otimes 1, 
  \q \ep_H^{\mrm{res}}=\ep_H \otimes 1, 
  \q S_H^{\mrm{res}}=S_H \otimes 1. 
\end{eqnarray*}

We define 
\begin{eqnarray*}
  X_{\pm, \AA}^{\mrm{res}}:=\sum_{j \in I, r \in \bbZ, m \in \bbN} 
  \AA (X_{j,r}^{\pm})^{(m)}, 
  \q X_{\pm}^{\mrm{res}}:=X_{\pm, \AA}^{\mrm{res}} \otimes 1
  =\sum_{j \in I, r \in \bbZ, m \in \bbN} \bbC (x_{j,r}^{\pm})^{(m)}.
\end{eqnarray*}
\begin{lem}
\label{lem DeH}
Let $i \in I$ and $r \in \bbN$. 
Modulo $\tUres \otimes \tUres X_{+}^{\mrm{res}}$,
\begin{eqnarray*}
  \De_H^{\mrm{res}}(x_{i,r}^{-})
  =x_{i,r}^{-} \otimes k_i + 1 \otimes x_{i,r}^{-}
  + \sum_{p=1}^{r-1} x_{i,r-p}^{-}  \otimes \psi_{i,p}^{+}.
\end{eqnarray*}
\end{lem}
Proof. 
It is enough to prove that 
for $i \in I$ and $r \in \bbZ$, 
modulo $\tUAres \otimes \tUAres X_{+, \AA}^{\mrm{res}}$,
\begin{eqnarray*}
  \De_H(X_{i,r}^{-})
  =X_{i,r}^{-} \otimes K_i + 1 \otimes X_{i,r}^{-}
  + \sum_{p=1}^{r-1} X_{i,r-p}^{-}  \otimes \Psi_{i,p}^{+}.
\end{eqnarray*}
From Proposition \ref{pro C02-2.2} (b), 
modulo $\tUq \otimes \tUq X_{+}$,
\begin{eqnarray*}
  \De_H(X_{i,r}^{-})
  =X_{i,r}^{-} \otimes K_i + 1 \otimes X_{i,r}^{-}
  + \sum_{p=1}^{r-1} X_{i,r-p}^{-}  \otimes \Psi_{i,p}^{+}.
\end{eqnarray*}
On the other hand, 
since $\tUAres$ is a Hopf subalgebra of $\tUres$, we have
\begin{eqnarray*}
  \De_H(X_{i,r}^{-}) \in \tUAres \otimes \tUAres.
\end{eqnarray*}
So it is enough to prove that
\begin{eqnarray*}
  (\tUq \otimes \tUq X_{+}) \cap (\tUAres \otimes \tUAres)
  =\tUAres \otimes \tUAres X_{+, \AA}^{\mrm{res}}. 
\end{eqnarray*}
This follows from (\ref{fac TD-GQLA}) and (\ref{TD-AQLA}). 
\qed
\subsection{Representation theory of $\Ures$ and $\tUres$}
We define 
\begin{eqnarray}
P_l:=\{\la=\sum_{i \in I} \la_i\La_i \in P \, | 
\, 0 \leq \la_i \leq l-1 \,\,  \rm{for any $i \in I$}\}. 
\label{def Pl}
\end{eqnarray}
For $\mu=\sum_{i \in I}\mu_i\La_i \in P$, 
there exist $\mu^{(0)}=\sum_{i \in I}\mu_i^{(0)}\La_i \in P_l$ 
and $\mu^{(1)}=\sum_{i \in I}\mu_i^{(1)}\La_i \in P$ 
such that $\mu=\mu^{(0)}+l\mu^{(1)}$. 
\begin{df} 
\label{def HWM-RQA}
Let $V$ be a representation of $\tUres$ (resp. $\Ures$). 

(i) Let $v \in V$. 
If $(x_{i,r}^{+})^{(m)}v=0$ 
for all $i \in I$, $r \in \bbZ$, and $m \in \bbN$ 
(resp. $e_i^{(m)}v=0$ for all $i \in I$ and $m \in \bbN$), 
we call $v$ a \it{pseudo-primitive vector} 
(resp. \it{primitive vector}) in $V$. 

(ii) For $\mu \in P$, we define 
\begin{eqnarray*}
  V_{\mu}:=\{v \in V \, 
  | \, k_iv=\vep^{\mu_i^{(0)}}v, \q
   \left[
  \begin{array}{r}
  k_i ; 0\\
  l \, \,
  \end{array}
  \right]v
  =\mu_i^{(1)}v  
  \q \rm{for all $i \in I$}\}. 
\end{eqnarray*}
If $V_{\mu} \neq 0$, we call $V_{\mu}$ a \it{weight space} of $V$. 
For $v \in V_{\mu}$, we call $v$ a weight vector with weight $\mu$ 
and define $\mrm{wt}(v):=\mu$.

(iii) For any $\bbC$-algebra homomorphism 
$\La: (\tUres)^0 \arr \bbC$, we define 
\begin{eqnarray*}
V_{\La}:=\{v \in V \, | \, uv=\La(u)v
\q \rm{for all $u \in (\tUres)^0$}\}. 
\end{eqnarray*}
If $V_{\La} \neq 0$, we call $V_{\La}$ a \it{pseudo-weight space} of $V$. 
For $v \in V_{\La}$, we call $v$ 
a pseudo-weight vector with pseudo-weight $\La$ 
and define $\mrm{pwt}(v):=\La$.

(iv) Let $\La: (\tUres)^0 \arr \bbC$ be a $\bbC$-algebra homomorphism 
and $\la$ be an element in $P_{+}$. 
If there exists a nonzero pseudo-primitive vector $v_{\La} \in V_{\La}$ 
(resp. primitive vector $v_{\la} \in V_{\la}$)
such that $V=\tUres v_{\La}$ (resp. $V=\Ures v_{\la}$), 
we call $V$ a \it{pseudo-highest weight representation} of $\tUres$ 
(resp. \it{highest-weight representation} of $\Ures$) 
with the \it{pseudo-highest weight $\La$} (resp. \it{highest weight $\la$})
generated by a \it{pseudo-highest weight vector $v_{\La}$} 
(resp. \it{highest-weight vector $v_{\la}$}). 
\end{df}

Let $V$ be a representation of $\tUres$ (resp. $\Ures$). 
We call $V$ of \it{type $\mbf{1}$} if 
$k_i^l=1$ on $V$ for any $i \in I$. 
In general, finite-dimensional irreducible representations 
of $\tUres$ (resp. $\Ures$) 
are classified into $2^n$ types according to 
$\{\sigma : Q \arr \{\pm 1\}; \rm{ group homomorphism}\}$. 
It is known that for any $\sigma : Q \arr \{\pm 1\}$, 
the category of finite-dimensional representations 
of $\tUres$ (resp. $\Ures$) of type $\sigma$ is essentially equivalent to 
the category of the finite-dimensional representations 
of $\tUres$ (resp. $\Ures$) of type $\mbf{1}$.
 
We define a set of polynomials $\bbC_0[t]$ whereby
\begin{eqnarray*}
\bbC_0[t]:=\{\pi(t) \in \bbC[t] \, | \, \pi(0)=1\}. 
\label{def RbbC0}
\end{eqnarray*}

For $\mbf{\pi}=(\pi_i(t))_{i \in I} \in (\bbC_0[t])^n$ 
(resp. $\la \in P_{+}$), 
let $v_{\mbf{\pi}}$ (resp. $v_{\la}$) be a 
pseudo-highest weight vector in $\tVq(\mbf{\pi})$ 
(resp. highest-weight vector in $\Vq(\la)$) 
(see \S 3.2).  
We define
\begin{eqnarray*}
&&\tVAres(\mbf{\pi}):=\tUAres v_{\mbf{\pi}}, \q
  \wt{W}_{\vep}^{\mrm{res}}(\mbf{\pi})
  :=\tVAres(\mbf{\pi}) \otimes_{\AA} \bbC_{\vep}, \\ 
&&\VAres(\la):=\UAres v_{\la}, \q
  W_{\vep}^{\mrm{res}}(\la):=\VAres(\la) \otimes_{\AA} \bbC_{\vep}. 
\end{eqnarray*}
We have
\begin{eqnarray*}
  \mrm{dim}_{\bbC}(\wt{W}_{\vep}^{\mrm{res}}(\mbf{\pi}))
  =\mrm{dim}_{\bbC(q)}(\tVq(\mbf{\pi})), 
  \q \mrm{dim}_{\bbC}(W_{\vep}^{\mrm{res}}(\la))
  =\mrm{dim}_{\bbC(q)}(\Vq(\la)),  
\end{eqnarray*}
(see \cite{CP94b}, Proposition 11.2.5).
For any $\tUres$-representation (resp. $\Ures$-representation) $V$, 
we have
\begin{eqnarray}
  (x_{i,r}^{\pm})^{(m)}V_{\mu} \subset V_{\mu \pm m\al_i}, 
  \q (\rm{resp.} \q e_i^{(m)} V_{\mu} \subset V_{\mu+m\al_i}, 
  \q f_i^{(m)} V_{\mu} \subset V_{\mu-m\al_i}), 
\label{fac Rweight space}
\end{eqnarray}
where $i \in I$, $r \in \bbZ$, $m \in \bbN$, and $\mu \in P$. 
Hence, by using (\ref{TD-RQLA}), 
we have the following proposition. 
\begin{pro}
\label{pro existence-HWR-RQLA}
For any $\mbf{\pi} \in (\bbC_0[t])^n$ (resp. $\la \in P_{+}$), 
$\wt{W}_{\vep}^{\mrm{res}}(\mbf{\pi})$ 
(resp. $W_{\vep}^{\mrm{res}}(\la)$) 
has a unique irreducible quotient 
$\tVres(\mbf{\pi})$ (resp. $\Vres(\la)$). 
\end{pro}
For any $\mbf{\pi}=(\pi_i(t))_{i \in I} \in (\bbC_0[t])^n$, 
there exists a unique $\bbC$-algebra homomorphism 
$\La_{\mbf{\pi}}^{\mrm{res}}: (\tUres)^0 \arr \bbC$ such that 
\begin{eqnarray*}
  \La_{\mbf{\pi}}^{\mrm{res}}(k_i^{\pm 1})
  =\vep^{\pm \mrm{deg}\pi_i(t)}, \q
  \sum_{m=1}^{\infty} \La_{\mbf{\pi}}^{\mrm{res}}(\frp_{i,\pm m})t^m 
  =\pi_i^{\pm}(t), 
  \label{def RLapi}
\end{eqnarray*}
where $\pi_i^{\pm}(t)$ be as in (\ref{def pi-}). 

For any pseudo-highest weight representation of $\tUres$ 
with the pseudo-highest weight $\La_{\mbf{\pi}}^{\mrm{res}}$, 
we simply call it a pseudo-highest representation of $\tUres$ 
with the pseudo-highest weight $\mbf{\pi}$. 
\begin{thm}[\cite{L89}]
\label{thm CT-RQA}
For $\la \in P_{+}$, 
$\Vres(\la)$ is a finite-dimensional irreducible 
highest-weight representation of $\Ures$ with the highest weight $\la$ 
of type $\mbf{1}$. 
Conversely, 
for any finite-dimensional irreducible $\Ures$-representation $V$ 
of type $\mbf{1}$, 
there exists a unique $\la \in P_{+}$ 
such that $V$ is isomorphic to $\Vres(\la)$ 
as a representation of $\Ures$.
\end{thm}
\begin{thm}[\cite{CP97}, \S 8]
\label{thm CT-RQLA}
For $\mbf{\pi} \in (\bbC_0[t])^n$, 
$\tVres(\mbf{\pi})$ is a finite-dimensional irreducible 
pseudo-highest weight representation of $\tUres$ 
with the pseudo-highest weight $\mbf{\pi}$ of type $\mbf{1}$. 
Conversely, 
for any finite-dimensional irreducible $\tUres$-representation $V$ 
of type $\mbf{1}$, 
there exists a unique $\mbf{\pi} \in (\bbC_0[t])^n$ 
such that $V$ is isomorphic to $\tVres(\mbf{\pi})$ 
as a representation of $\tUres$.
\end{thm}
\begin{pro}[\cite{CP97}, Proposition 7.4]
\label{pro tensor highest-RQA}
Let $V$ (resp. $V^{'}$) be a pseudo-highest weight representation 
of $\tUres$ with the pseudo-highest weight 
$\mbf{\pi}$ (resp. $\mbf{\pi}^{'}$) 
generated by a pseudo-highest weight vector 
$v_{\mbf{\pi}}$ (resp. $v_{\mbf{\pi}^{'}}^{'}$). 
Then $v_{\mbf{\pi}} \otimes v_{\mbf{\pi}^{'}}^{'}$ is 
a pseudo-primitive vector with 
$\mrm{pwt}(v_{\mbf{\pi}} \otimes v_{\mbf{\pi}^{'}}^{'})
=\La_{\mbf{\pi} \mbf{\pi}^{'}}^{\mrm{res}}$.
\end{pro}
From Proposition 8.3 in \cite{CP97} and 
Proposition \ref{pro tensor highest-RQA}, 
we obtain the following corollary.
\begin{cor}
\label{cor tensor highest-RQA}
Let $\mbf{\pi}, \mbf{\pi}^{'} \in (\bbC_0[t])^n$ and 
let $v_{\mbf{\pi}}$ (resp. $v_{\mbf{\pi}^{'}}^{'}$) 
be a pseudo-highest weight vector in 
$\tVres(\mbf{\pi})$ (resp. $\tVres(\mbf{\pi}^{'})$). 
$\tVres(\mbf{\pi}\mbf{\pi^{'}})$ 
is isomorphic to a quotient of the $\tUres$-subrepresentation 
of $\tVres(\mbf{\pi}) \otimes \tVres(\mbf{\pi}^{'})$ 
generated by $v_{\mbf{\pi}} \otimes v_{\mbf{\pi}^{'}}^{'}$. 
In particular, 
if $\tVres(\mbf{\pi}) \otimes \tVres(\mbf{\pi}^{'})$ is irreducible, 
$\tVres(\mbf{\pi}) \otimes \tVres(\mbf{\pi}^{'})$ 
is isomorphic to $\tVres(\mbf{\pi}\mbf{\pi}^{'})$. 
\end{cor}
\subsection{The fundamental representations of $\Ures$ and $\tUres$}
For $\xi \in I$ and $\bf{a} \in \bbC^{\times}$, 
let $\mbf{\pi}_{\xi}^{\mbf{a}}
=(\pi_{\xi,j}^{\mbf{a}}(t))_{j \in I} \in (\bbC_0[t])^n$ 
be as in (\ref{def FP}). 
We call $\tVres(\mbf{\pi}_{\xi}^{\mbf{a}})$ 
(resp. $\Vres(\La_{\xi})$) 
a \it{fundamental representation} of $\tUres$ (resp. $\Ures$).
We have
\begin{eqnarray}
  \VAres(\La_{\xi})=\bigoplus_{\mu \in \WW\La_{\xi}} \AA z_{\mu},
\label{fac VAres}
\end{eqnarray}
and $W_{\vep}^{\mrm{res}}(\La_{\xi})
=\VAres(\La_{\xi}) \otimes_{\AA} \bbC_{\vep}$ 
is irreducible as a representation of $\Ures$. 
So we identify $W_{\vep}^{\mrm{res}}(\La_{\xi})$ 
with $\Vres(\La_{\xi})$. 
For $z \in \VAres(\La_{\xi})$, we set 
\begin{eqnarray*}
 \ol{z}:=z \otimes 1 \in \Vres(\La_{\xi}).
\end{eqnarray*}
From (\ref{fac VAres}) and Lemma \ref{lem A-basis}, we have
\begin{eqnarray}
  \Vres(\La_{\xi})=\bigoplus_{\mu \in \WW\La_{\xi}} \bbC \ol{z}_{\mu}.
\label{fac Vres}
\end{eqnarray}
For $\xi \in I$ and $\mbf{a} \in \bbC^{\times}$, 
let $\VAres(\La_{\xi})_{\mbf{a}}$ 
be the $\tUAres$-subrepresentation of $\Vq(\La_{\xi})_{\mbf{a}}$ 
generated by $z_{\La_{\xi}}$. 
We define 
$\Vres(\La_{\xi})_{\mbf{a}}
:=\VAres(\La_{\xi})_{\mbf{a}} \otimes_{\AA} \bbC_{\vep}$.
Then, as representations of $\tUres$, we have 
\begin{eqnarray}
  \tVres(\mbf{\pi}_{\xi}^{\mbf{a}}) \cong \Vres(\La_{\xi})_{\mbf{a}}. 
\label{fac dim}
\end{eqnarray} 

From (\ref{pro Hs1}), for $i \in I$, we have
\begin{eqnarray}
  h_{i,1}\ol{z}_{\La_{\xi}}
  =\La_{\mbf{\pi}_{\xi}^{\mbf{a}}}^{\mrm{res}}(h_{i,1})
  \ol{z}_{\La_{\xi}}
  =\mbf{a}\vep^{-1} \de_{i, \xi}\ol{z}_{\La_{\xi}} 
  \q \rm{in} \q \Vres(\La_{\xi})_{\mbf{a}}.
\label{fac aez}
\end{eqnarray}
Moreover, from Proposition \ref{pro extremal}, 
we obtain the following proposition.
\begin{pro} 
\label{pro Rextremal}
Let $\xi \in I$ and $\mbf{a} \in \bbC^{\times}$. 
For $i, j \in I$ such that $i \leq j$, 
let $\om_{i,j}$ be as in (\ref{def om}) 
and let $\ol{z_{\La_{\xi}}(\om_{i,j})}$ be the 
extremal vector in $\Vres(\La_{\xi})_{\mbf{a}}$. 

(a) We have $\ol{z_{\La_{\xi}}(\om_{i,j})}
    =\ol{z_{\om_{i,j}\La_{\xi}}}$.

(b) We have
\begin{eqnarray*}
&& h_{i-1,1}\ol{z_{\La_{\xi}}(\om_{i,j})}
  =\mbf{a}\vep^{2j-i-\xi}\de(j-i+2 \leq \xi \leq j)
  \ol{z_{\La_{\xi}}(\om_{i,j})}, \q \rm{if $i \neq 1$},\\
&& h_{j+1,1}\ol{z_{\La_{\xi}}(\om_{1,j})}
  =\mbf{a}\vep^{j-\xi}\de(1 \leq \xi \leq j+1)
  \ol{z_{\La_{\xi}}(\om_{1,j})}.
\end{eqnarray*}
\end{pro}
From (\ref{fac Om}), 
we have a $\bbC$-algebra involution 
$\Om^{\mrm{res}}: \tUres \arr \tUres$ such that 
\begin{eqnarray*}
 \Om^{\mrm{res}}((x_{i,r}^{\pm})^{(m)})=-(x_{i,-r}^{\mp})^{(m)}, 
\q \Om^{\mrm{res}}(h_{i,s})=-h_{i,-s}, 
\q \Om^{\mrm{res}}(\psi_{i,r}^{\pm})=\psi_{i,-r}^{\mp}, 
\q \Om^{\mrm{res}}(k_i^{\pm 1})=k_i^{\mp 1},
\end{eqnarray*}
for any $i \in I$, $r \in \bbZ$, $s \in \bbZ^{\times}$, 
and $m \in \bbN$.
From Proposition \ref{pro DFR+IFR-GQLA}, 
we have the following proposition.
\begin{pro} 
\label{pro R-DFR+IFR-GQLA}
Let $\xi \in I$ and $\mbf{a} \in \bbC^{\times}$. 

(a) There exists an integer $c \in \bbZ$ depending only $\sl_{n+1}$ 
    such that, as a representation of $\tUres$, 
    $\Vres(\La_{\xi})_{\mbf{a}}^{*}$ is isomorphic to 
    $\Vres(\La_{n-\xi+1})_{\vep^{c} \mbf{a}}$.    

(b) There exists a nonzero complex number $\kappa \in \bbC^{\times}$ 
    depending only $\sl_{n+1}$ such that, as a representation of $\tUres$, 
    $\Vres(\La_{\xi})_{\mbf{a}}^{\Om^{\mrm{res}}}$ is isomorphic to 
    $\Vres(\La_{n-\xi+1})_{\vep^2\kappa \mbf{a}^{-1}}$. 
    In particular, $\ol{z_{\om_{1,n}\La_{\xi}}}$ 
    is a pseudo-highest weight vector in 
    $\Vres(\La_{\xi})_{\mbf{a}}^{\Om^{\mrm{res}}}$.
\end{pro}

By using Proposition \ref{pro tensor highest-RQA} and 
Proposition \ref{pro R-DFR+IFR-GQLA} (b), 
in a similar way to the proof of Lemma \ref{lem HWR+IR=highest}, 
we obtain the following lemma.
\begin{lem}
\label{lem R-HWR+IR=highest}
Let $\xi \in I$, $\mbf{a} \in \bbC^{\times}$, 
and $\mbf{\pi} \in (\bbC_0[t])^n$. 
Let $V$ be a pseudo-highest weight representation of $\tUq$
with the pseudo-highest weight $\mbf{\pi}$ 
and let $v_{\pi}$ be a pseudo-highest weight vector in $V$. 
We assume $\ol{z_{\om_{1,n}\La_{\xi}}} \otimes v 
\in \tUres(\ol{z_{\La_{\xi}}} \otimes v)$. 
Then $\Vres(\La_{\xi})_{\mbf{a}} \otimes V$ 
is a pseudo-highest weight representation of $\tUres$ 
with the pseudo-highest weight 
$\mbf{\pi}_{\xi}^{\mbf{a}}\mbf{\pi}$.
\end{lem}
\subsection{Irreducibility: the $\tUres$ case}
For $r \in \bbZ$, $s \in \bbZ^{\times}$, and $m \in \bbN$, 
we denote $(x_{1,r}^{\pm})^{(m)}$ (resp. $h_{1,s}$, $k_1^{\pm 1}$) 
in $\Ures(\wt{\sl}_2)$ by $(x_r^{\pm})^{(m)}$ (resp. $h_s$, $k^{\pm 1}$) 
and $e_1$ (resp. $f_1$, $k_1^{\pm 1}$) in $\Ures(\sl_2)$ 
by $e$ (resp. $f$, $k^{\pm 1}$). 
For $i \in I$, 
let $(\tUres)^{(i)}$ (resp. $(\Ures)^{(i)}$) 
be the $\bbC$-subalgebra of $\tUres$ (resp. $\Ures$) 
generated by $\{(x_{i,r}^{\pm})^{(m)}, k_i^{\pm 1} \, |
\, r \in \bbZ, m \in \bbN\}$ 
(resp. $\{e_i^{(m)}, f_i^{(m)}, k_i^{\pm 1} \, | \, m \in \bbN\}$). 
There exist $\bbC$-algebra homomorphisms 
$\wt{\iota}^{\mrm{res}}: \Ures(\wt{\sl}_2) \arr (\tUres)^{(i)}$ 
and $\iota^{\mrm{res}}: \Ures(\sl_2) \arr (\Ures)^{(i)}$ such that 
\begin{eqnarray*}
&& \wt{\iota}^{\mrm{res}}((x_r^{\pm})^{(m)})=(x_{i,r}^{\pm})^{(m)}, \q 
\wt{\iota}^{\mrm{res}}(h_s)=h_{i,s}, \q 
\wt{\iota}^{\mrm{res}}(k^{\pm 1})=k_i^{\pm 1}, \\
&& \iota^{\mrm{res}}(e^{(m)})=e_i^{(m)}, 
\q \iota^{\mrm{res}}(f^{(m)})=f_i^{(m)}, 
\q \iota^{\mrm{res}}(k^{\pm 1})=k_i^{\pm 1},
\end{eqnarray*}
(see \S 5.1).
Hence, for any $(\tUres)^{(i)}$-representation 
(resp. $(\Ures)^{(i)}$-representation) $V$, 
we can regard $V$ as a $\Ures(\wt{\sl}_2)$-representation 
(resp. $\Ures(\sl_2)$-representation). 

In a similar way to the proof of Lemma \ref{lem TP-FR-GQLA}, 
we can prove the following lemma.
\begin{lem}
\label{lem R-TP-FR-GQLA}
Let $\xi \in I$ and $\mbf{a} \in \bbC^{\times}$. 
For any $i, j \in I$ such that $i \leq j$,
let $\ol{z_{\om_{i, j}\La_{\xi}}}$ be the extremal vector in 
$\Vres(\La_{\xi})_{\mbf{a}}$.
As representations of $\Ures(\wt{\sl}_2)$, 
\begin{eqnarray*}
&&(\tUres)^{(i-1)} \ol{z_{\om_{i,j}\La_{\xi}}} \cong
  \begin{cases}
    \Vres(1)_{\mbf{a}\vep^{2j-\xi-i+1}},
    & \rm{if} \q j-i+2 \leq \xi \leq j, \\
    \Vres(0)_{\mbf{a}}, & otherwise, 
  \end{cases}
  \q \rm{if} \q i \neq 1,\\
&&(\tUres)^{(j+1)} \ol{z_{\om_{1,j}\La_{\xi}}} \cong
  \begin{cases}
    \Vres(1)_{\mbf{a} \vep^{j-\xi+1}},
    & \rm{if} \q 1 \leq \xi \leq j+1, \\
    \Vres(0)_{\mbf{a}}, & otherwise. 
  \end{cases}
\end{eqnarray*}
\end{lem}
\begin{thm}
\label{thm TP-FR-RQLA}  
Let $m \in \bbN$, $\xi_1, \cdots, \xi_m \in I$, 
and $\mbf{a}_1, \cdots, \mbf{a}_m \in \bbC^{\times}$. 
We assume that for any $1 \leq k < k^{'} \leq m$ 
and $\mrm{max}(\xi_k, \xi_{k^{'}}) \leq t 
\leq \mrm{min}(\xi_k+\xi_{k^{'}}-1, n)$, 
\begin{eqnarray*}
\frac{\mbf{a}_{k^{'}}}{\mbf{a}_k}
\neq \vep^{2t-\xi_k-\xi_{k^{'}}+2}. 
\end{eqnarray*}
Then $\Vres(\La_{\xi_1})_{\mbf{a}_1} \otimes 
\cdots \otimes \Vres(\La_{\xi_m})_{\mbf{a}_m}$ 
is a pseudo-highest weight representation of $\tUres$
with the pseudo-highest weight 
$\mbf{\pi}^{\mbf{a}_1}_{\xi_1} \cdots \mbf{\pi}^{\mbf{a}_m}_{\xi_m}$ 
generated by a pseudo-highest weight vector 
$\ol{z_{\La_{\xi_1}}} \otimes \cdots \otimes \ol{z_{\La_{\xi_m}}}$.
\end{thm}

Proof. 
We can prove this theorem in a similar way to the proof of 
Theorem \ref{thm TP-FR-GQLA}. 
 
From Proposition \ref{pro tensor highest-RQA}, 
it is enough to prove 
\begin{eqnarray*}
  \Vres(\La_{\xi_1})_{\mbf{a}_1} \otimes 
  \cdots \otimes \Vres(\La_{\xi_m})_{\mbf{a}_m}
  =\tUres (\ol{z_{\La_{\xi_1}}} \otimes \cdots 
  \otimes \ol{z_{\La_{\xi_m}}}).
\end{eqnarray*}
We shall prove this claim by the induction on $m$. 
If $m=1$, we have nothing to prove. 
So we assume $m>1$ and the case of $(m-1)$ holds. 
We set 
\begin{eqnarray*}
  V^{'}:=\Vres(\La_{\xi_2})_{\mbf{a}_2} \otimes 
  \cdots \otimes \Vres(\La_{\xi_m})_{\mbf{a}_m}, 
  \q z^{'}:=\ol{z_{\La_{\xi_2}}} \otimes \cdots 
  \otimes \ol{z_{\La_{\xi_m}}}.
\end{eqnarray*}
From Proposition \ref{pro tensor highest-RQA} 
and the assumption of the induction on $m$, 
$V^{'}$ is a pseudo-highest weight representation of $\tUres$
with the pseudo-highest weight 
$\mbf{\pi}^{\mbf{a}_2}_{\xi_2} \cdots \mbf{\pi}^{\mbf{a}_m}_{\xi_m}$ 
generated by a pseudo-highest weight vector $z^{'}$. 
Hence, from Lemma \ref{lem R-HWR+IR=highest}, 
it is enough to prove that 
\begin{eqnarray*}
  \ol{z_{\om_{1,n}\La_{\xi_1}}} \otimes z^{'} 
  \in \tUres(\ol{z_{\La_{\xi_1}}} \otimes z^{'}).
  \end{eqnarray*}
We shall prove that
\begin{eqnarray}
  \ol{z_{\om_{i,j}\La_{\xi_1}}} \otimes z^{'} 
  \in \tUres(\ol{z_{\La_{\xi_1}}} \otimes z^{'}),
\label{thm R-TP-FR-GQLA claim}
\end{eqnarray} 
for any $i, j \in I$ such that $i \leq j$. 
We define a total order in $I^{\leq}$ as (\ref{def TO2}).
We shall prove (\ref{thm R-TP-FR-GQLA claim}) by the induction on $(i,j)$. 
If $(i,j)=(1,0)$, we have nothing to prove. 
So we assume that the case of $(i,j)$ holds. 
We also assume $i \neq 1$. 
We can prove the case of $i=1$ similarly.
If $\xi_1<j-i+2$ or $\xi_1 > j$, we have
\begin{eqnarray*}
  \ol{z_{\om_{i-1,j}\La_{\xi_1}}} \otimes z^{'}
  =\ol{z_{\om_{i,j}\La_{\xi_1}}} \otimes z^{'}
  \in \tUres(\ol{z_{\La_{\xi_1}}} \otimes z^{'}).  
\end{eqnarray*}
So we assume $j-i+2 \leq \xi_1 \leq j$.
From Lemma \ref{lem DeH}, for $r \in \bbN$, we have
\begin{eqnarray*}
&&\De_{H}^{\mrm{res}}(x_{i-1,r}^{-})
  (\ol{z_{\om_{i,j}\La_{\xi_1}}} \otimes z^{'})
  -\ol{z_{\om_{i,j}\La_{\xi_1}}} \otimes 
  (\De_{H}^{\mrm{res}}(x_{i-1,r}^{-})z^{'}) \no \\ 
&&=(\De_{H}^{\mrm{res}}(x_{i-1,r-k}^{-})\ol{z_{\om_{i,j}\La_{\xi_1}}})
  \otimes (\De_{H}^{\mrm{res}}(k_{i-1})z^{'})
  +\sum_{k=1}^{r-1}(x_{i-1,r-k}^{-}\ol{z_{\om_{i,j}\La_{\xi_1}}})
  \otimes (\De_{H}^{\mrm{res}}(\psi_{i-1,k}^{+})z^{'}).
\end{eqnarray*}
(see (\ref{thm TP-FR-GQLA 4})).
By using Lemma \ref{lem R-TP-FR-GQLA} and Lemma \ref{lem detA}, 
we obtain
\begin{eqnarray*}
  (\prod_{k \in M} (\mbf{a}_k-\mbf{a}_1 \vep^{2j-\xi_1-\xi_k+2}))
  (\prod_{k, k^{'} \in M, k<k^{'}} 
  (\mbf{a}_{k^{'}}-\mbf{a}_k \vep^2))
  (\ol{z_{\om_{i-1,j}\La_{\xi_1}}} \otimes z^{'})
  \in \tUres (\ol{z_{\La_{\xi_1}}} \otimes z^{'}),
\label{thm R-TP-FR-GQLA 11}
\end{eqnarray*}
where $M$ be as in (\ref{thm TP-FR-GQLA 12}). 
Therefore, from the assumption of this theorem, we obtain 
\begin{eqnarray*}
  \ol{z_{\om_{i-1,j}\La_{\xi_1}}} \otimes z^{'}
  \in \tUres (\ol{z_{\La_{\xi_1}}} \otimes z^{'}).
\end{eqnarray*}
\qed \\
By using Theorem \ref{thm TP-FR-RQLA} and 
Proposition \ref{pro R-DFR+IFR-GQLA}, 
we obtain the following corollary 
(see Corollary \ref{cor TP-FR-GQLA}). 
\begin{cor} 
\label{cor R-TP-FR-GQLA}
Let $m \in \bbN$, $\xi_1, \cdots, \xi_m \in I$, 
and $\mbf{a}_1, \cdots, \mbf{a}_m \in \bbC^{\times}$. 
We assume that for any $1 \leq k \neq k^{'} \leq m$ 
and $\mrm{max}(\xi_k, \xi_{k^{'}}) \leq t 
\leq \mrm{min}(\xi_k+\xi_{k^{'}}-1, n)$, 
\begin{eqnarray*}
\frac{\mbf{a}_{k^{'}}}{\mbf{a}_k}
\neq \vep^{\pm (2t-\xi_k-\xi_{k^{'}}+2)}. 
\end{eqnarray*}
Then $\Vres(\La_{\xi_1})_{\mbf{a}_1} \otimes 
\cdots \otimes \Vres(\La_{\xi_m})_{\mbf{a}_m}$ 
is an irreducible representation of $\tUres$.
\end{cor}
\subsection{Reducibility: the $\tUres$ case}
\begin{pro} 
\label{pro R-reducibility}
Let $\xi, \ze \in I$ and $\mbf{a}, \mbf{b} \in \bbC^{\times}$. 
If there exists a 
$1 \leq t \leq \mrm{min}(\xi, \ze, n+1-\xi, n+1-\ze)$ such that 
$\mbf{b}=\vep^{2t+|\xi-\ze|}\mbf{a}$ or $\vep^{-(2t+|\xi-\ze|)}\mbf{a}$, 
then $\Vres(\La_{\xi})_{\mbf{a}} \otimes \Vres(\La_{\ze})_{\mbf{b}}$ 
is reducible as a representation of $\tUres$. 
\end{pro}

Proof. 
Let $q$ be an indeterminate 
and let 
\begin{eqnarray*}
  (\mbf{b}_q, \mbf{b}_{\vep})=
  (q^{2t+|\xi-\ze|}\mbf{a}, \vep^{2t+|\xi-\ze|}\mbf{a})  \q \rm{or} 
  \q (q^{-(2t+|\xi-\ze|)}\mbf{a}, \vep^{-(2t+|\xi-\ze|)}\mbf{a}). 
\end{eqnarray*}
We assume that  
$\Vres(\La_{\xi})_{\mbf{a}} \otimes \Vres(\La_{\ze})_{\mbf{b}_{\vep}}$ 
is irreducible as a representation of $\tUres$. 
Since $\Vres(\La_{\xi})_{\mbf{a}} \cong
\tVres(\mbf{\pi}_{\xi}^{\mbf{a}})$ 
(resp. $\Vres(\La_{\ze})_{\mbf{b}_{\vep}} \cong 
\tVres(\mbf{\pi}_{\ze}^{\mbf{b}_{\vep}})$), 
from Corollary \ref{cor tensor highest-RQA}, we obtain
\begin{eqnarray*}
  \Vres(\La_{\xi})_{\mbf{a}} \otimes \Vres(\La_{\ze})_{\mbf{b}_{\vep}} 
  \cong \tVres (\mbf{\pi}_{\xi}^{\mbf{a}}\mbf{\pi}_{\ze}^{\mbf{b}_{\vep}}). 
\end{eqnarray*}
Hence, from (\ref{fac Vres}) and (\ref{fac dim}), we have
\begin{eqnarray*}
  \dim_{\bbC}(\tVres (\mbf{\pi}_{\xi}^{\mbf{a}}
  \mbf{\pi}_{\ze}^{\mbf{b}_{\vep}}))
  =\dim_{\bbC} (\Vres(\La_{\xi})_{\mbf{a}}) 
  \times \dim_{\bbC}(\Vres(\La_{\ze})_{\mbf{b}_{\vep}}) 
  =\dim_{\bbC(q)} (\Vq(\La_{\xi})_{\mbf{a}}) 
  \times \dim_{\bbC(q)}(\Vq(\La_{\ze})_{\mbf{b}_q}).
\end{eqnarray*} 
On the other hand, by the definition of 
$\tVres(\mbf{\pi}_{\xi}^{\mbf{a}}\mbf{\pi}_{\ze}^{\mbf{b}_{\vep}})$, 
we have 
\begin{eqnarray*}
  \dim_{\bbC}(\tVres(\mbf{\pi}_{\xi}^{\mbf{a}}
  \mbf{\pi}_{\ze}^{\mbf{b}_{\vep}}))   \leq 
  \dim_{\bbC(q)} \tVq(\mbf{\pi}_{\xi}^{\mbf{a}}
  \mbf{\pi}_{\ze}^{\mbf{b}_q}),    
\end{eqnarray*}
(see \S 6.4). 
Thus, we have 
\begin{eqnarray*}
  \dim_{\bbC(q)} (\Vq(\La_{\xi})_{\mbf{a}}) 
  \times \dim_{\bbC(q)}(\Vq(\La_{\ze})_{\mbf{b}_q}) \leq 
  \dim_{\bbC(q)} (\tVq(\mbf{\pi}_{\xi}^{\mbf{a}}
  \mbf{\pi}_{\ze}^{\mbf{b}_q})). 
\end{eqnarray*}
Hence, from Corollary \ref{cor tensor highest}, 
\begin{eqnarray*}
\Vq(\La_{\xi})_{\mbf{a}}  \otimes \Vq(\La_{\ze})_{\mbf{b}_q}
\cong \tVq(\mbf{\pi}_{\xi}^{\mbf{a}}\mbf{\pi}_{\ze}^{\mbf{b}_q}).  
\end{eqnarray*}
In particular, 
$\Vq(\La_{\xi})_{\mbf{a}}  \otimes \Vq(\La_{\ze})_{\mbf{b}_q}$ 
is irreducible as a representation of $\tUq$. 
However, from Proposition \ref{pro Reducibility-FR-GQLA}, 
$\Vq(\La_{\xi})_{\mbf{a}}  \otimes \Vq(\La_{\ze})_{\mbf{b}_q}$ 
is reducible. 
This is absurd. 
Therefore 
$\Vres(\La_{\xi})_{\mbf{a}} \otimes \Vres(\La_{\ze})_{\mbf{b}_{\vep}}$ 
is reducible as a representation of $\tUres$. 
\qed
\subsection{Main theorem: the $\tUres$ case}
\begin{thm} 
\label{thm Main theorem-RQLA}
Let $m \in \bbN$, $\xi_1, \cdots, \xi_m \in I$, 
and $\mbf{a}_1, \cdots, \mbf{a}_m \in \bbC^{\times}$. 
The following conditions (a) and (b) are equivalent. 

(a) 
$\tVres(\pi_{\xi_1}^{\mbf{a}_1}) \otimes 
\cdots \otimes \tVres(\pi_{\xi_m}^{\mbf{a}_m})$ 
is an irreducible representation of $\tUres$.

(b) 
For any $1 \leq k \neq k^{'} \leq m$ 
and $1 \leq t \leq 
\mrm{min}(\xi_k, \xi_{k^{'}}, n+1-\xi_k, n+1-\xi_{k^{'}})$, 
\begin{eqnarray*}
  \frac{\mbf{a}_{k^{'}}}{\mbf{a}_k}
  \neq \vep^{\pm (2t+|\xi_k-\xi_{k^{'}}|)}. 
\end{eqnarray*}
\end{thm}

Proof. 
This theorem follows from 
Corollary \ref{cor R-TP-FR-GQLA}
and Proposition \ref{pro R-reducibility}
(see Theorem \ref{thm Main theorem-GQLA}).
\qed 
\section{Tensor product of the fundamental representations for 
the small quantum loop algebras}
\setcounter{equation}{0}
\renewcommand{\theequation}{\thesection.\arabic{equation}}
\subsection{The tensor product theorems}
Let $P_l$ be as in (\ref{def Pl}). 
For $\la \in P_{+}$, 
let $\la^{(0)} \in P_l$ and $\la^{(1)} \in P_{+}$ 
be as in \S 6.4. 
\begin{thm}[\cite{L89}, Theorem 7.4]
\label{thm TPT-RQA}
For $\la \in P_{+}$, 
$\Vres(\la)$ is isomorphic to 
$\Vres(\la^{(0)}) \otimes \Vres(l\la^{(1)})$ 
as a representation of $\Ures$. 
\end{thm}
For $\pi(t) \in \bbC_0[t]$, we call $\pi(t)$ \it{$l$-acyclic} 
if it is not divisible by $(1-ct^l)$ for any $c \in \bbC^{\times}$ 
(see \cite{FM}, \S 2.6).
We define 
\begin{eqnarray}
&&\bbC_l[t]:=\{\pi(t) \in \bbC_0[t] \, | \, 
  \rm{$\pi(t)$ is $l$-acyclic}\}, 
  \label{def Cl}\\
&&\bbC[t^l]:=\{\pi(t) \in \bbC_0[t] \, | \, 
  \rm{there exists a polynomial $\pi^{'}(t) \in \bbC_0[t]$ 
  such that $\pi(t)=\pi^{'}(t^l)$}\}.
\label{def Cl2}
\end{eqnarray}
For $\mbf{\pi} \in (\bbC_0[t])^n$, 
there exist unique 
$\mbf{\pi}^{(0)}=(\pi_i^{(0)}(t))_{i \in I} \in (\bbC_l[t])^n$ 
and $\mbf{\pi}^{(1)}=(\pi_i^{(1)}(t))_{i \in I} \in (\bbC_0[t^l])^n$ 
such that $\pi_i(t)=\pi_i^{(0)}(t)\pi_i^{(1)}(t)$ for any $i \in I$. 
\begin{thm}[\cite{CP97}, Theorem 9.1]
\label{thm TPT-RQLA}
For $\mbf{\pi} \in (\bbC_0[t])^n$, 
$\tVres(\mbf{\pi})$ is isomorphic to 
$\tVres(\mbf{\pi}^{(0)}) \otimes \tVres(\mbf{\pi}^{(1)})$ 
as a representation of $\tUres$.
\end{thm}
\subsection{The Frobenius homomorphisms 
and the construction of $\Vres(l\la)$ and $\tVres(\mbf{\pi}^{(1)})$}

Let $\wt{U}:=U(\wt{\sl}_{n+1})$ (resp. $U:=U(\sl_{n+1})$) 
be the universal enveloping algebra of 
$\wt{\sl}_{n+1}$ (resp. $\sl_{n+1}$), 
that is, $\wt{U}$ (resp. $U$) is 
an associative algebra over $\bbC$ 
generated by 
$\{\bar{e}_i, \bar{f}_i, \bar{h}_i \, | \, 
i \in \wt{I} \, \, (\rm{resp.} \, \, i \in I)\}$ 
with the defining relations:
\begin{eqnarray*}
&& \bar{h}_i\bar{h}_j=\bar{h}_j\bar{h}_i, 
\q \bar{h}_i\bar{e}_j-\bar{e}_j\bar{h}_i
=\fra_{i,j}\bar{e}_j, 
\q \bar{h}_i\bar{f}_j-\bar{f}_j\bar{h}_i
=-\fra_{i,j}\bar{f}_j, 
\q \bar{e}_i\bar{f}_j-\bar{f}_j\bar{e}_i
=\de_{i,j}\bar{h}_i, \\
&& \sum_{m=0}^{1-\fra_{i,j}}
(-1)^m
\frac{(1-\fra_{i,j})!}{m!(1-\fra_{i,j}-m)!}
 \bar{e}_i^{m}\bar{e}_j\bar{e}_i^{1-\fra_{i,j}-m}
=\sum_{m=0}^{1-\fra_{i,j}} 
(-1)^m
\frac{(1-\fra_{i,j})!}{m!(1-\fra_{i,j}-m)!}
\bar{f}_i^{m}\bar{f}_j\bar{f}_i^{1-\fra _{i,j}-m}
=0 \q i \neq j.
\end{eqnarray*}
Then, from \cite{CP97}, \S 1, we have the following theorem 
(see also \cite{CP94b}, Theorem 9.3.12 and \S 11.2B).
\begin{thm}[\cite{CP97}, \S 1] 
\label{thm FM}
There exist a $\bbC$-algebra homomorphism
$\wt{\mrm{Fr}}_{\vep}: \tUres \arr \wt{U}$ 
such that 
\begin{eqnarray*}
&&\wt{\mrm{Fr}}_{\vep}(e_i^{(m)})=
\begin{cases}
\frac{\bar{e}_i^{m/l}}{(m/l)!}, & \rm{if $l$ divides $m$}, \\
0, & \rm{otherwise},
\end{cases}
\q \wt{\mrm{Fr}}_{\vep}(f_i^{(m)})=
\begin{cases}
\frac{\bar{f}_i^{m/l}}{(m/l)!}, & \rm{if $l$ divides $m$}, \\
0, & \rm{otherwise},
\end{cases} \\
&&\wt{\mrm{Fr}}_{\vep}(k_i)=1, \q
\wt{\mrm{Fr}}_{\vep}([k_i;l])=\bar{h}_i,
\end{eqnarray*}
for any $i \in I$ and $m \in \bbN$.
\end{thm}
For any $\wt{U}$-representation $V$, 
we can regard $V$ as a $\tUres$-representation 
by using $\wt{\mrm{Fr}}_{\vep}$ 
and denote it by $\wt{\mrm{Fr}}_{\vep}^{*}(V)$. 
Similarly, for any $U$-representation $V$, 
we can regard $V$ as a $\Ures$-representation by using 
$\mrm{Fr}_{\vep}:=\wt{\mrm{Fr}}_{\vep}|_{\Ures}: \Ures \arr U$ 
and denote it by $\mrm{Fr}_{\vep}^{*}(V)$. 
For $\mbf{\pi} \in (\bbC_0[t])^n$ (resp. $\la \in P_{+}$), 
let $\wt{V}(\mbf{\pi})$ (resp. $V(\la)$) 
be the finite-dimensional irreducible representation of 
$\wt{U}$ (resp. $U$) with the pseudo-highest weight $\mbf{\pi}$ 
(resp. highest-weight $\la$) (see \cite{CP97}, \S 2).
\begin{thm}[\cite{L89}, \S 7 and \cite{CP94b}, Proposition 11.2.11]
\label{thm FH-RQA}
For $\la \in P_{+}$, 
$\Vres(l\la)$ is isomorphic to $\mrm{Fr}_{\vep}^{*}(V(\la))$ 
as a representation of $\Ures$. 
\end{thm}
\begin{thm}[\cite{CP97}, Theorem 9.3]
\label{thm FH-RQLA}
For $\mbf{\pi}=(\pi_i(t))_{i \in I}, 
\mbf{\pi}^{'}=(\pi_i^{'}(t))_{i \in I} \in (\bbC_0[t])^n$ 
such that $\pi_i(t)=\pi_i^{'}(t^l)$ for any $i \in I$, 
$\tVres(\mbf{\pi})$ is isomorphic to 
$\wt{\mrm{Fr}}_{\vep}^{*}(V(\pi^{'}))$ 
as a representation of $\tUres$. 
\end{thm}
\subsection{Definition and the representation theory 
of the small quantum algebras}
\begin{df} 
\label{def SQA}
Let $\tUfin$ (resp. $\Ufin$) be the $\bbC$-subalgebra of 
$\tUres$ (resp. $\Ures$) generated by 
$\{e_i, f_i, k_i \, | \, i \in \wt{I} \, (\mrm{resp.} \, i \in I)\}$. 
We call $\tUfin$ a \it{small quantum loop algebra} 
(resp. \it{small quantum algebra}).
\end{df}
For any $\tUfin$-representation (resp. $\Ufin$-representation) $V$, 
we call $V$ of \it{type $\mbf{1}$} 
if $k_i=1$ on $V$ for all $i \in I$. 
For $\pi \in (\bbC_0[t])^n$ (resp. $\la \in P_{+}$), 
we regard $\tVres(\pi)$ (resp. $\Vres(\la)$) 
as a representation of $\tUfin$ (resp. $\Ufin$) 
and denote it by $\tVfin(\pi)$ (resp. $\Vfin(\la)$).
\begin{thm}[\cite{L89}, Proposition 7.1 and 
\cite{CP94b}, Proposition 11.2.10]
\label{thm CT-SQA}
For $\la \in P_l$, 
$\Vfin(\la)$ is irreducible as a representation of $\Ufin$. 
Moreover,  
for any finite-dimensional irreducible $\Ufin$-representation $V$ 
of type $\mbf{1}$, 
there exists a unique $\la \in P_l$ 
such that $V$ is isomorphic to $\Vfin(\la)$ 
as a representation of $\Ufin$.
\end{thm}
\begin{thm}[\cite{CP97}, Theorem 9.2 and \cite{FM}, Thorem 2.6]
\label{thm CT-SQLA}
For $\mbf{\pi} \in (\bbC_l[t])^n$, 
$\tVfin(\mbf{\pi})$ is irreducible as a representation of $\tUfin$. 
Moreover, 
for any finite-dimensional irreducible $\tUfin$-representation $V$ 
of type $\mbf{1}$, 
there exists a unique $\mbf{\pi} \in (\bbC_l[t])^n$ 
such that $V$ is isomorphic to $\tVfin(\mbf{\pi})$ 
as a representation of $\tUfin$.
\end{thm}
\begin{rem}
\label{rem CT}
From Theorem \ref{thm TPT-RQLA} and Theorem \ref{thm FH-RQLA} 
(resp. Theorem \ref{thm TPT-RQA} and Theorem \ref{thm FH-RQA}), 
in order to understand the finite-dimensional irreducible 
representations of $\tUres$ (resp. $\Ures$), 
we may consider the one of $\tUfin$ (resp. $\Ufin$). 
\end{rem}
\subsection{Main theorem: the $\tUfin$ case}
\begin{thm} 
\label{thm Main theorem-SQLA}
Let $m \in \bbN$, $\xi_1, \cdots, \xi_m \in I$, 
and $\mbf{a}_1, \cdots, \mbf{a}_m \in \bbC^{\times}$. 
The following conditions (a) and (b) are equivalent. 

(a) 
$\tVfin(\pi_{\xi_1}^{\mbf{a}_1}) \otimes 
\cdots \otimes \tVfin(\pi_{\xi_m}^{\mbf{a}_m})$ 
is an irreducible representation of $\tUfin$.

(b) 
For any $1 \leq k \neq k^{'} \leq m$ 
and $1 \leq t \leq 
\mrm{min}(\xi_k, \xi_{k^{'}}, n+1-\xi_k, n+1-\xi_{k^{'}})$, 
\begin{eqnarray*}
  \frac{\mbf{a}_{k^{'}}}{\mbf{a}_k}
  \neq \vep^{\pm (2t+|\xi_k-\xi_{k^{'}}|)}. 
\end{eqnarray*}
\end{thm}
Proof. 
If (b) does not hold, then (a) also does not hold 
from Theorem \ref{thm Main theorem-RQLA}. 
So we assume that (b) holds. 
From Theorem \ref{thm Main theorem-RQLA} and 
Corollary \ref{cor tensor highest-RQA}, 
as representations of $\tUres$,
\begin{eqnarray*}
  \tVfin(\pi_{\xi_1}^{\mbf{a}_1}) \otimes 
  \cdots \otimes \tVfin(\pi_{\xi_m}^{\mbf{a}_m}) 
  =\tVres(\pi_{\xi_1}^{\mbf{a}_1}) \otimes 
  \cdots \otimes \tVres(\pi_{\xi_m}^{\mbf{a}_m})
  \cong \tVres(\mbf{\pi}_{\xi_1}^{\mbf{a}} 
  \cdots \mbf{\pi}_{\xi_m}^{\mbf{a}_m}). 
\end{eqnarray*} 
Hence, from Theorem \ref{thm CT-SQLA}, 
it is enough to prove 
$\mbf{\pi}_{\xi_1}^{\mbf{a}}  \cdots \mbf{\pi}_{\xi_m}^{\mbf{a}_m} 
\in (\bbC_l[t])^n$. 

There exist $\pi_i(t) \in \bbC_0[t]$ ($i \in I$) such that 
$\mbf{\pi}_{\xi_1}^{\mbf{a}}  \cdots \mbf{\pi}_{\xi_m}^{\mbf{a}_m} 
=(\pi_i(t))_{i \in I}$. 
If there exists an index $i \in I$ 
such that $\pi_i(t) \not\in \bbC_l[t]$, 
there exists a nonzero complex number $c$ such that 
$(1-ct)(1-c\vep t) \cdots (1-c\vep^{l-1}t)$ divides $\pi_i(t)$. 
Then there exist $1 \leq i_1, \cdots, i_t \leq m$ such that 
\begin{eqnarray*}
  \xi_{i_1}= \cdots =\xi_{i_t}, \q
  \mbf{a}_{i_1}=c, \q \mbf{a}_{i_2}=c\vep,  \q \cdots, \q 
  \mbf{a}_{i_t}=c\vep^{l-1}. 
\end{eqnarray*} 
On the other hand, since (b) holds, 
\begin{eqnarray*}
  \frac{\mbf{a}_{i_s}}{\mbf{a}_{i_r}} \neq \vep^{\pm 2}, 
\end{eqnarray*}
for any $1 \leq r \neq s \leq t$.
This is absurd.
\qed\\

\bf{Acknowledgements:} 
I would like to thank Masaharu Kaneda, Hyohe Miyachi, 
Toshiki Nakashima, and Atsushi Nakayashiki 
for their helpful discussions.


\end{document}